\documentclass{article}%
\usepackage{amsmath}
\usepackage{amssymb}
\usepackage{amsmath}
\usepackage{amsfonts}
\usepackage{graphicx}
\usepackage{stmaryrd}
\DeclareMathOperator*{\esssup}{ess\,sup}

\usepackage[english]{babel}
\usepackage{dsfont}
\usepackage{url}
\usepackage{authblk} 

\newtheorem{theorem}{Theorem}

\newtheorem{corollary}[theorem]{Corollary}

\newtheorem{definition}[theorem]{Definition}

\newtheorem{lemma}[theorem]{Lemma}

\newtheorem{proposition}[theorem]{Proposition}
\newtheorem{remark}[theorem]{Remark}

\newenvironment{proof}[1][Proof]{\noindent\textbf{#1.} }{\ \rule{0.5em}{0.5em}}
\usepackage{vmargin}
\setmarginsrb{3.5cm}{4cm}{3.5cm}{4cm}{0cm}{0cm}{0cm}{0.5cm}
\begin{document}
	\title{$\mathbb{L}^p$-solution of generalized BSDEs in a general filtration with stochastic monotone coefficients}
	%
	
\author[1]{Badr ELMANSOURI\thanks{Corresponding author\\	This research was supported by the National Center for Scientific and Technical Research (CNRST), Morocco.}} 
	\author[1]{Mohamed EL OTMANI} 
	\affil[1]{Laboratory of Analysis, Mathematics, and Applications (LAMA), 
		Faculty of Sciences Agadir, 
		Ibnou Zohr University, Agadir, Morocco} 
	
	\affil[ ]{\texttt{badr.elmansouri@edu.uiz.ac.ma}} 
	\affil[ ]{\texttt{m.elotmani@uiz.ac.ma}} 
	\date{}
	\maketitle
	\begin{abstract}
		We study multidimensional generalized backward stochastic differential equations (GBSDEs) within a general filtration that supports a Brownian motion under weak assumptions on the associated data. We establish the existence and uniqueness of solutions in $\mathbb{L}^p$ for $p \in (1,2]$. Our results apply to generators that are stochastic monotone in the $y$-variable, stochastic Lipschitz in the $z$-variable, and satisfy a general stochastic linear growth condition.
	\end{abstract}
	
	
	
	\section{Introduction}
	The theory of backward stochastic differential equations (BSDEs) has been thoroughly studied and shown to have a wide range of applications in various mathematical domains, including partial differential equations (PDEs) \cite{Pardoux1999}, stochastic control and differential games \cite{hamadene1995backward,Hamadene1997}, mathematical finance \cite{el1997backward,yong2006completeness}, and other related fields. Bismut \cite{bismut1973conjugate} originally introduced the concept as the adjoint equations related to stochastic Pontryagin maximum principles in stochastic control theory. Pardoux and Peng \cite{pardoux1990adapted} were the first to study the general case of nonlinear multidimensional BSDEs. Roughly speaking, for a finite horizon time $T \in (0,+\infty)$, a solution of such equations, associated with a terminal value $\xi$ and a generator (coefficient or driver) $f(\omega,t,y,z)$, is a pair of stochastic processes $(Y_t,Z_t)_{t \leq T}$ satisfying:
	\begin{equation}\label{basic BSDE}
		Y_t=\xi+\int_{t}^{T}f(s,Y_s,Z_s)ds-\int_{t}^{T}Z_s dW_s,~t \in [0,T],
	\end{equation}
	where $W=(W_t)_{t \leq T}$ is a standard Brownian motion, and the solution process $(Y_t,Z_t)_{t \leq T}$ is adapted to the natural filtration of $W$. Under a uniform Lipschitz condition on the driver $f$ and a square integrability condition on $\xi$ and the process $(f(\omega,t,0,0))_{t \leq T}$, the authors \cite{pardoux1990adapted} demonstrated the existence and uniqueness of a solution. 
	
	Following this work, many researchers have aimed to weaken the uniform Lipschitz continuity constraint on the generator to address more interesting problems. In this context, significant research has been conducted on the existence, uniqueness, and comparison theorems for $\mathbb{L}^2$-solutions of the BSDE (\ref{basic BSDE}) with square-integrable parameters and a weaker condition than the Lipschitz one considered in \cite{pardoux1990adapted}; see, e.g., \cite{bahlali2001backward,fan2010uniqueness,hamadene2003multidimensional}, among others. However, in some practical applications, even when considering an appropriate condition on the generator $f$ weaker than the Lipschitz one, the terminal condition $\xi$ and the driver process $(f(\omega,t,0,0))_{t \leq T}$ of the BSDE (\ref{basic BSDE}) are not necessarily assumed to be square-integrable. Consequently, considering BSDEs with $\mathbb{L}^p$-integrable data and $\mathbb{L}^p$-solutions for $p \geq 1$ has attracted significant interest over the last decade. Briand et al. \cite{briand2003lp} demonstrated the existence and uniqueness of $\mathbb{L}^p$-solutions for $p \in (1,2)$ of the BSDE (\ref{basic BSDE}) when the generator $f$ is monotonic in $y$ and Lipschitz continuous in $z$. For additional relevant works, see \cite{chen2010p,eddahbi,fan2012p,fan2010class,yao2017lp} and the references therein.

	Using a new class of BSDEs that involves the integral with respect to a continuous increasing process interpreted as the local time of a diffusion process on the boundary, Pardoux and Zhang \cite{pardoux1998generalized} provided a probabilistic representation for a solution of a system of parabolic and elliptic semi-linear PDEs with Neumann boundary conditions. This new type of BSDE is called Generalized BSDEs (GBSDEs). A solution of such equations is a pair of adapted processes $(Y_t,Z_t)_{t \leq T}$ that satisfies the following equation:
	\begin{equation}\label{GBSDE}
		Y_t=\xi+\int_{t}^{T}f(s,Y_s,Z_s)ds+\int_{t}^{T}g(s,Y_s)d\kappa_s-\int_{t}^{T}Z_s dW_s,~t \in [0,T].
	\end{equation}
	Here:
	\begin{itemize}
		\item $\xi$ is an $\mathbb{R}^d$-valued $\mathcal{F}_T$-measurable random variable,
		\item $f : \Omega \times [0,T] \times \mathbb{R}^d \times \mathbb{R}^{d \times k} \rightarrow \mathbb{R}^d$ is an $\mathcal{F} \otimes \mathcal{B}([0,T])\otimes \mathcal{B}(\mathbb{R}^d)\otimes \mathcal{B}(\mathbb{R}^{d \times k})$-measurable random function such that for any $(y,z) \in \mathbb{R}^d \times \mathbb{R}^{d \times k}$, the process $(\omega,t) \mapsto f(\omega,t,y,z)$ is progressively measurable, 
		\item $g : \Omega \times [0,T] \times \mathbb{R}^d \rightarrow \mathbb{R}^d$ is an $\mathcal{F} \otimes \mathcal{B}([0,T])\otimes \mathcal{B}(\mathbb{R}^d)$-measurable random function such that for any $y \in \mathbb{R}^d$, the process $(\omega,t) \mapsto g(\omega,t,y)$ is progressively measurable, 
		\item $\kappa=(\kappa_t)_{t \leq T}$ is an $\mathbb{R}$-valued adapted, continuous, increasing process on $[0,T]$.
	\end{itemize} 
	Under a monotonicity assumption on the drivers $f$ and $g$, and appropriate $\mathbb{L}^2$-integrability conditions on the data, the authors established the existence and uniqueness of a solution using a convolution approximation. Extending this framework, Pardoux \cite{pardoux1997generalized} addressed the discontinuous case by incorporating a jump term into (\ref{GBSDE}), represented by an independent Poisson random measure. More recently, Elmansouri and El Otmani \cite{ElmansouriElOtmani,elmansourielotmani2024} demonstrated existence and uniqueness results for GBSDEs in a general filtration under similar or more general assumptions compared to those in \cite{pardoux1997generalized,pardoux1998generalized}.
	
	However, all the aforementioned works, concerning BSDE (\ref{basic BSDE}) or GBSDE (\ref{GBSDE}), deal with square or $\mathbb{L}^p$-integrable parameters and different weak conditions on the drivers only in a Brownian framework. In such a case, it is well known that the predictable representation property holds for every local martingale (see, e.g., Theorem 43 in \cite[p. 186]{protter2005stochastic}). However, this is no longer valid for more general filtrations (see Section III.4 in \cite{jacod2003limit}), and the description of a solution must include an extra orthogonal martingale term to $W$. More precisely, for a filtration $\mathbb{F}:=(\mathcal{F}_t)_{t \leq T}$ carrying (or supporting) the Brownian motion $W$, every right-continuous with left limits (RCLL) local martingale $(\mathcal{N}_t)_{t \leq T}$ can be represented as follows (see, e.g., Lemma 4.24 in \cite[p. 185]{jacod2003limit}):
	$$
	\mathcal{N}_t=\int_{0}^{t}Z_s dW_s+M_t,\quad t \in [0,T],
	$$
	for some predictable process $(Z_t)_{t \leq T}$ such that $\int_{0}^{T}\|Z_s\|^2 ds<+\infty$ a.s. and an RCLL local martingale $(M_t)_{t \leq T}$ such that $\left[M,W\right]=0$ a.s., where $[M,W]$ denotes the co-variation process of $M$ and $W$. For RCLL martingales, this method was developed in the groundbreaking works of Carbone et al. \cite{carbone2008backward} and El Karoui and Huang \cite{ElKarouiHuang1997} for classical BSDEs, and later by Elmansouri and El Otmani \cite{elmansourielotmani2024} for GBSDEs in a more general framework. In the same filtration context but within the $\mathbb{L}^p$-setup, Kruse and Popier \cite{kruse2016bsdes,kruse2017lp} studied the $\mathbb{L}^p$-solution ($p>1$) problem for BSDEs in a general filtration supporting a Brownian motion and an independent Poisson random measure. The authors proved the existence and uniqueness of a solution with a driver $f$ that is monotone with respect to $y$ and uniformly Lipschitz with respect to $z$.
	
	Compared to the existing literature and to the best of our knowledge, $\mathbb{L}^p$-solutions for $p \in (1,2)$ of GBSDEs in the general filtration $\mathbb{F}$ have not been widely investigated. More precisely, we consider the following multidimensional GBSDE:
	\begin{equation}\label{basic GBSDE}
		Y_t=\xi+\int_{t}^{T}f(s,Y_s,Z_s)ds+\int_{t}^{T}g(s,Y_s)d\kappa_s-\int_{t}^{T}Z_s dW_s-\int_{t}^{T}dM_s,~t \in [0,T].
	\end{equation}

	Motivated by the works mentioned above, it is naturally interesting to investigate the existence and uniqueness results for $\mathbb{L}^p$-solutions ($p \in (1,2)$) of multidimensional GBSDEs (\ref{basic GBSDE}) under suitable, more general conditions on the data. To this end, under a stochastic monotonicity condition on $f$ and $g$ with respect to the $y$-variable, a stochastic Lipschitz condition on $f$ with respect to the $z$-variable, a general stochastic linear growth condition, and an appropriate $\mathbb{L}^p$-integrability condition on the data, we aim to establish a general existence and uniqueness result for $\mathbb{L}^p$-solutions of the multidimensional GBSDEs (\ref{basic GBSDE}) for $p \in (1,2)$ and also for $p = 2$, as we could not find a result for the latter case under our general assumptions.
	
	In the sequel, we present several publications that studied the uniqueness and existence of $\mathbb{L}^p$-solutions for BSDEs (\ref{basic BSDE}) with time-varying or stochastic monotonic conditions on the coefficient $f$, since results for the GBSDE (\ref{basic GBSDE}) have not been examined previously. Xiao et al. \cite{xiao2015lp} consider the case where the driver $g$ satisfies a time-varying monotonicity condition in $y$, meaning that there exists a deterministic integrable function $[0,T] \ni t \mapsto \alpha_t \in \mathbb{R}_+$ such that for each $y, y' \in \mathbb{R}^d$ and $z \in \mathbb{R}^{d \times k}$,
	$$
	(y - y')\left(f(t, y, z) - f(t, y', z)\right) \leq \alpha_t \left| y - y' \right|^2, \quad d\mathbb{P} \otimes dt \text{-a.e.},
	$$
	and a time-varying Lipschitz continuity condition on $z$, meaning that there exists a deterministic square-integrable function $[0,T] \ni t \mapsto \eta_t \in \mathbb{R}_+$ such that, for any $z, z' \in \mathbb{R}^{d \times k}$,
	$$
	\left| f(t, y, z) - f(t, y, z') \right| \leq \eta_t \left\| z - z' \right\|, \quad d\mathbb{P} \otimes dt \text{-a.e.}.
	$$
	Under these conditions and appropriate $\mathbb{L}^p$-integrability conditions on the data $\xi$ and the process $(f(t, 0, 0))_{t \leq T}$, the authors in \cite{xiao2015lp} prove the existence and uniqueness of $\mathbb{L}^p$-solutions for $p > 1$ using the method of convolution and weak convergence. Pardoux and R\^{a}\c{s}canu
	 \cite{pardoux2014stochastic} also study existence and uniqueness results for $\mathbb{L}^p$-solutions ($p > 1$) for multidimensional BSDEs of the form (\ref{basic BSDE}) and (\ref{GBSDE}) (see \cite[Chapter 5]{pardoux2014stochastic}) under different growth conditions, including the case where $(\alpha_t)_{t \leq T}$, $(\eta_t)_{t \leq T}$ are deterministic functions and where $\alpha$ takes values in $\mathbb{R}$ as stochastic processes. They also provide the connection with semi-linear PDEs and parabolic variational inequalities with a mixed nonlinear multi-valued Neumann–Dirichlet boundary condition. Very recently, Li et al. \cite{li2024weighted} provided existence and uniqueness of $\mathbb{L}^p$-solutions for the BSDE (\ref{basic BSDE}) under a stochastic monotonicity condition on the driver $f$ with respect to $(y, z)$, as a direct extension of the work by Li et al. \cite{li2024random} in the $\mathbb{L}^2$ case.
	
	In this paper, under the above-mentioned stochastic monotonicity and Lipschitz conditions on the drivers $f$ and $g$, we establish a general existence and uniqueness result for $\mathbb{L}^p$-solutions ($p \in (1,2]$) of multidimensional GBSDEs in a general filtration carrying a Brownian motion. The first part of this paper is devoted to establishing essential a priori $\mathbb{L}^p$-estimates of the solutions to the GBSDE (\ref{basic GBSDE}) for $p \in (1,2)$. It is worth noting that our state process $(Y_t)_{t \leq T}$ is not necessarily continuous but only RCLL, which introduces additional challenges in our work. Specifically, compared to (\ref{basic BSDE}) or (\ref{GBSDE}), our GBSDE (\ref{basic GBSDE}) includes a jump term represented by the orthogonal martingale $M$, complicating the proof since the bracket process involving the quadratic jumps of $M$ (or the state process $Y$) must be carefully handled. Afterward, using these results, we study the existence and uniqueness of $\mathbb{L}^p$-solutions for $p \in (1,2]$ when the generator $f$ and $g$ are stochastic monotonic with respect to $y$, and $f$ is stochastic Lipschitz with respect to $z$, along with a general stochastic linear growth condition in $y$. We first treat the case $p=2$, breaking the proof into two main steps. The first considers the case when $f$ depends only on $y$, where we apply the Yosida approximation method. Then, using this result and a fixed-point theorem, we extend the case to a general driver $f$. From the $\mathbb{L}^2$ case, we derive the existence and uniqueness results for $\mathbb{L}^p$ solutions ($p \in (1,2)$) using an appropriate sequence of GBSDEs of the form (\ref{basic GBSDE}).
	
	It is worth mentioning that, by using the results presented in this paper and following the arguments in \cite[Section 6]{kruse2016bsdes}, \cite[Theorem 4.1]{Pardoux1999} or \cite[Theorem 53.2]{pardoux1998generalized}, our work can be extended to the case where the terminal time $T$ is replaced by a stopping time $\tau$ for the general filtration $\mathbb{F}$, which need not be bounded. Additionally, the case $p \in (2,+\infty)$ can also be treated within our framework without additional complexity compared to the $p \in (1,2)$ case, due to the fact that for $p > 2$, the function $\mathbb{R}^d \ni x \mapsto |x|^p$ is sufficiently smooth, allowing direct application of It\^o's formula and other classical arguments (see \cite{eddahbi} for a related study). The case $p \in (1,2)$ is less regular and requires other representation formulas.

	The rest of this paper is organized as follows: Section \ref{sec2} introduces some notations, definitions, and results used in the paper. Section \ref{sec3} establishes some important a priori $\mathbb{L}^p$-estimates ($p \in (1,2)$) for solutions of the GBSDE (\ref{basic GBSDE}). In Section \ref{sec4}, we prove the existence and uniqueness result for the $\mathbb{L}^p$-solutions for $p \in (1,2]$.

	\section{Preliminaries}
	\label{sec2}
	Let $T > 0$ be a fixed deterministic horizon time, and let $\left(\Omega, \mathcal{F}, \mathbb{P}\right)$ be a complete probability space equipped with a filtration $\mathbb{F} := (\mathcal{F}_t)_{t \leq T}$, carrying a $k$-dimensional Brownian motion $(W_t)_{t \leq T}$. The filtration $\mathbb{F}$ is assumed to satisfy the usual conditions of right-continuity and completeness. The initial $\sigma$-field $\mathcal{F}_0$ is assumed to be trivial, and $\mathcal{F}_T = \mathcal{F}$. Unless explicitly stated, all stochastic processes are considered on the time interval $[0,T]$, and the measurability properties of stochastic processes (such as adaptedness, predictability, progressive measurability) are taken with respect to $\mathbb{F}$.
	
	The bracket process (or quadratic variation) of any given $\mathbb{R}^d$-valued RCLL local martingale $M$ is defined by $\left[ M \right]$. The notation $\left[ M \right]^c$ denotes the continuous part of the quadratic variation $\left[ M \right]$. Specifically, the bracket process $\left[ M \right]$ of $M$ is defined for every $t \in [0,T]$ as follows:
	$$
	\left[ M \right]_t = \sum_{i=1}^{d} \left[ M^i \right]_t,
	$$
	where $M^i$ is the $i$th component of the vector $M = \left(M^1, M^2, \cdots, M^d\right)$.\\ 
	For another given RCLL local martingale $N$, $[M, N]$ denotes the quadratic covariation matrix process of $M$ and $N$, defined for every $t \in [0,T]$ as follows:
	$$
	[M, N]_t = \left\{\left[M^i, N^j\right]_t; 1 \leq i, j \leq d \right\}.
	$$
	
	The Euclidean norm of a vector $y \in \mathbb{R}^d$ is defined by $
	\left| y \right|^2 = \sum_{i=1}^{d} \left| y_i \right|^2,$ and for a given matrix $z \in \mathbb{R}^{k \times d}$, we set $
	\|z\|^2 = \mbox{Trace}(z z^\ast),$ where $z^\ast$ denotes the transpose of $z$.
	
	For an RCLL process $(\mathcal{X}_t)_{t \leq T}$, $\mathcal{X}_{t-} := \lim\limits_{s \nearrow t} \mathcal{X}_s$ denotes the left limit of $\mathcal{X}$ at $t$. We set $\mathcal{X}_{0-} = \mathcal{X}_0$ by convention. The process $\mathcal{X}_{-} = (\mathcal{X}_{t-})_{t \in [0,T]}$ is called the left-limited process, and $\Delta \mathcal{X} = \mathcal{X} - \mathcal{X}_{-}$ is the jump process associated with $\mathcal{X}$. More precisely, for any $t \in [0,T]$, we have $\Delta \mathcal{X}_t = \mathcal{X}_t - \mathcal{X}_{t-}$, which is the jump of $\mathcal{X}$ at time $t$.
	
	For an adapted process with finite variation $\mathcal{V} = (\mathcal{V}_t)_{t \leq T}$, we denote by $\|\mathcal{V}\| = (\|\mathcal{V}\|_t)_{t \leq T}$ the total variation process on $[0,T]$.
	
	To simplify the notation, we omit any dependence on $\omega$ of a given process or random function. By convention, all brackets and stochastic integrals are assumed to be zero at time zero.
	

	Let $\beta, \mu\geq 0$, $p>1$ and $(a_t)_{t \leq T}$ a progressively measurable positive process. We consider the increasing continuous process $A_t:=\int_{0}^{t}a_s^2 ds$ for $t \in [0,T]$. 
	
	To define the $\mathbb{L}^p$-solution of our GBSDE (\ref{basic GBSDE}) for $p>1$, we need to introduce the following spaces:
	\begin{itemize}
		\item $\mathcal{S}^p_{\beta}$: The space of $\mathbb{R}^d$-valued and $\mathcal{F}_{t}$-adapted RCLL processes $(Y_{t})_{t\leq T}$ such that
		$$
		\|Y\|_{\mathcal{S}^{p}_{\beta,\mu}}=\left(\mathbb{E}\left[\sup\limits_{t \in [0,T]}e^{\frac{p}{2}\beta A_t+\frac{p}{2}\mu \kappa_t}|Y_{t}|^{p}\right]\right)^{\frac{1}{p}} < +\infty.
		$$
		
		\item[$\bullet$] $\mathcal{S}^{p,A}_{\beta,\mu}$ is the space of $\mathbb{R}^d$-valued and $\mathcal{F}_{t}$-adapted RCLL processes $(Y_{t})_{t\leq T}$ such that
		$$
		\|Y\|_{\mathcal{S}^{p,A}_{\beta,\mu}}=\left(\mathbb{E}\left[\int_{0}^{T}e^{\frac{p}{2}\beta A_s+\frac{p}{2}\mu \kappa_s} |Y_{s}|^{p}dA_s\right]\right)^\frac{1}{p} < +\infty.
		$$
		
		\item[$\bullet$] $\mathcal{S}^{p,A}_{\beta,\mu}$ is the space of $\mathbb{R}^d$-valued and $\mathcal{F}_{t}$-adapted RCLL processes $(Y_{t})_{t\leq T}$ such that
		$$
		\|Y\|_{\mathcal{S}^{p,\kappa}_{\beta,\mu}}=\left(\mathbb{E}\left[\int_{0}^{T}e^{\frac{p}{2}\beta A_s+\frac{p}{2}\mu \kappa_s} |Y_{s}|^{p}d\kappa_s\right]\right)^\frac{1}{p} < +\infty.
		$$
		
		\item[$\bullet$] $\mathcal{H}^{p}_{\beta,\mu}$ is the space of $\mathbb{R}^{d\times k}$-valued and predictable processes $(Z_{t})_{t\leq T}$ such that
		$$\|Z\|_{\mathcal{H}^{p}_{\beta,\mu}}=\left(\mathbb{E}\left[\left(\int_{0}^{T}e^{\beta A_s+\mu \kappa_s}\|Z_{s}\|^{2}ds\right)^\frac{p}{2}\right]\right)^\frac{1}{p} < +\infty.$$
		
		\item[$\bullet$] $\mathcal{M}^{p}_{\beta,\mu}$ is the space of all martingales $(M_t)_{t \leq T}$ orthogonal to $W$ such that
		$$\|M\|_{\mathcal{M}^{p}_{\beta,\mu}}=\left(\mathbb{E}\left[\left(\int_{0}^{T}e^{\beta A_s+\mu \kappa_s} d\big[M\big]_s\right)^\frac{p}{2}\right]\right)^\frac{1}{p} < +\infty.$$
		
		Finally, we define 
		\item[$\bullet$] $\mathfrak{B}^p_{\beta,\mu}:=\mathcal{S}^{p}_{\beta,\mu}\cap \mathcal{S}^{p,A}_{\beta,\mu}\cap \mathcal{S}^{p,\kappa}_{\beta,\mu}$ endowed with the norm
		$\|Y\|_{\mathfrak{B}^{p}_{\beta,\mu}= \|Y\|_{\mathcal{S}^{p}_{\beta,\mu}}^p}+\|Y\|_{\mathcal{S}^{p,A}_{\beta,\mu}}^p+\|Y\|_{\mathcal{S}^{p,\kappa}_{\beta,\mu}}^p.$
		
		\item $\mathcal{E}^p_{\beta,\mu}:=\mathfrak{B}^p_{\beta,\mu} \times \mathcal{H}^p_{\beta,\mu}\times \mathcal{M}^p_{\beta,\mu}$ endowed with the norm $ \|(Y,Z,M)\|^p_{\mathcal{E}^p_{\beta,\mu}}:=\|Y\|_{\mathfrak{B}^{p}_{\beta,\mu}}^p+\|Z\|^p_{\mathcal{H}^{p}_{\beta,\mu}}+\|M\|^p_{\mathcal{M}^{p}_{\beta,\mu}}$. 
	\end{itemize}

	The following two remarks provide a simple statement that will be useful throughout the paper when dealing with $\mathbb{L}^p$-estimates for $p \in (1,2)$.
	\begin{remark}\label{Good rmq too}
		Let $M$ and $N$ be two $\mathbb{R}^d$-valued RCLL local martingales. Using the polarization identity for the bracket process $\left[M - N\right]$ (see, e.g., \cite[p. 52]{jacod2003limit}) along with its definition given above, we get 
		\begin{equation*}
			\begin{split}
				\left[M - N\right] = \sum_{i=1}^d \left[M^i - N^i\right] 
				&= \sum_{i=1}^d \left\{ \left[M^i\right] - 2\left[M^i, N^i\right] + \left[N^i\right] \right\} \\
				&= \left[M\right] + \left[N\right] - 2 \sum_{i=1}^d \left[M^i, N^i\right].
			\end{split}
		\end{equation*}
		Now assume that $M, N \in \mathcal{M}^p_{\beta, \mu}$ for $\beta, \mu > 0$. By the Kunita-Watanabe inequality (see Theorem 25 in \cite[p. 69]{protter2005stochastic}), for any $\tau \in \mathcal{T}_{0,T}$ and all $i \in \{1, \dots, d\}$, we have $\mathbb{P}$-a.s.
		\begin{equation*}
			\begin{split}
				&-2 \int_0^\tau e^{\beta A_s + \mu {\kappa}_s} d \left[M^i, N^i\right]_s\\ &\leq 2 \left( \int_0^\tau e^{\beta A_s +\mu {\kappa}_s} d \left[M^i\right]_s \right)^{\frac{1}{2}} \left( \int_0^\tau e^{\beta A_s + \mu {\kappa}_s} d \left[N^i\right]_s \right)^{\frac{1}{2}} \\
				&\leq \left( \int_0^\tau e^{\beta A_s + \mu {\kappa}_s} d \left[M^i\right]_s + \int_0^\tau e^{\beta A_s + \mu {\kappa}_s} d \left[N^i\right]_s \right).
			\end{split}
		\end{equation*}
		Therefore, after summing over all $i \in \{1, \dots, d\}$, we derive
		\begin{equation*}
			\begin{split}
				&\mathbb{E}\left[\left(\int_{0}^{\tau} e^{\beta A_s + \mu {\kappa}_s} d \left[M-N\right]_s\right)^{\frac{p}{2}}\right]\\
				& \leq 2^p	\left(\mathbb{E} \left[\left(\int_{0}^{\tau} e^{\beta A_s + \mu {\kappa}_s} d \left[M\right]_s\right)^{\frac{p}{2}}\right]+	\mathbb{E}\left[\left(\int_{0}^{\tau} e^{\beta A_s + \mu {\kappa}_s} d \left[N\right]_s\right)^{\frac{p}{2}}\right]\right)<+\infty.
			\end{split}
		\end{equation*}
	\end{remark}
	
	Let $p > 1$. Throughout the rest of this paper, the following conditions on the terminal value $\xi$ and the generator $f$ are denoted by \textbf{(H-M)$_\mathbf{p}$}.
	\paragraph{Conditions on the data $(\xi,f,g,\kappa)$}
	For any $\beta, \mu >0$, we assume that
	\begin{itemize}
		\item[\textsc{(H1)}] The terminal condition $\xi$ satisfies
		$$
		\mathbb{E}\left[e^{\frac{p}{2}\beta A_T+\frac{p}{2}\mu \kappa_T}\left|\xi\right|^p\right]<+\infty.
		$$
		\item[\textsc{(H2)}] There exists an $\mathbb{F}$-progressively measurable processes $\alpha : \Omega \times [0,T] \rightarrow \mathbb{R}$ such that for all $t \in [0,T]$, $y,y' \in \mathbb{R}^d$, $z \in \mathbb{R}^{d \times k}$, $d\mathbb{P}\otimes dt$-a.e.,
		$$
		\left(y-y'\right)\left(f(t,y,z)-f(t,y',z)\right) \leq \alpha_t \left|y-y'\right|^2.
		$$
		
		\item[\textsc{(H3)}] There exists an $\mathbb{F}$-progressively measurable processes $\theta : \Omega \times [0,T] \rightarrow \mathbb{R}_-^\ast$ such that for all $t \in [0,T]$, $y,y' \in \mathbb{R}^d$, $z \in \mathbb{R}^{d \times k}$, $d\mathbb{P}\otimes dt$-a.e.,
		$$
		\left(y-y'\right)\left(g(t,y)-g(t,y')\right) \leq \theta_t \left|y-y'\right|^2.
		$$
		
		\item[\textsc{(H4)}] There exists an $\mathbb{F}$-progressively measurable processes $\eta : \Omega \times [0,T] \rightarrow \mathbb{R}^\ast_{+}$ and  such that for all $t \in [0,T]$, $y, \in \mathbb{R}^d$, $z,z' \in \mathbb{R}^{d \times k}$, $d\mathbb{P}\otimes dt$-a.e.,
		$$
		\left| f(t,y,z)-f(t,y,z')\right|  \leq \eta_t \left\|z-z'\right\|.
		$$
		
		\item[\textsc{(H5)}] There exist three progressively measurable processes $\varphi, \psi : \Omega \times \left[0,T\right] \rightarrow [1,+\infty)$, $\phi : \Omega \times \left[0,T\right] \rightarrow (0,+\infty)$, $\zeta : \Omega \times \left[0,T\right] \rightarrow (0,\Gamma]$ for some $\Gamma>0$ such that $\left|f(t,y,0)\right| \leq \varphi_t+\phi_t |y|$, $\left|g(t,y)\right| \leq \psi_t+\zeta_t |y|$ and 
		$$
		\mathbb{E}{\int _{0}^{T}} e^{\beta A_s+\mu \kappa_s} \left(\left|\varphi_s\right|^pds+\left|\psi_s\right|^pd\kappa_s\right) <+\infty ,
		$$
		
		\item[\textsc{(H6)}] We set $a^2_s=\phi^2_s+\eta^2_s+\delta^2_s$ and we assume that there exists a constant $\epsilon>0$ such that $a^2_s \geq \epsilon$ for any $s \in [0,T]$.  
		
		\item[\textsc{(H7)}] For all $(t, z) \in [0,T] \times \mathbb{R}^{d \times k} $, the mappings $y \mapsto f(t, y, z)$ and  $y \mapsto g(t, y)$ are continuous, $d\mathbb{P} \otimes dt$-a.e. and $d\mathbb{P} \otimes d\kappa_t$-a.e., respectively.	
	\end{itemize}
	\begin{remark}
		Throughout this paper, $\mathfrak{c}$ denote a positive constant that may change from one line to another. Additionally, the notation $\mathfrak{c}_\gamma$ is used to emphasize the dependence of the constant $\mathfrak{c}$ on a specific set of parameters $\gamma$.
	\end{remark}

	To derive the optimal constants in the a priori estimates of the solutions, we use the following remark:
	\begin{remark}\label{rmq essential}
		Let $(Y,Z,M)$ be a solution of the GBSDE (\ref{basic GBSDE}) associated with $(\xi, f,g,\kappa)$. Let $\varepsilon \geq 0$, $\varrho \in \mathbb{R}$ and $(\lambda^{(\alpha),\varepsilon,\varrho}_t)_{t \leq T}$ a progressively measurable increasing continuous  process defined by ${\lambda}^{(\alpha),\varepsilon,\varrho}_t = \exp\left\{\int_{0}^{t}\alpha_s ds+\varepsilon A_t + \varrho \kappa_t\right\}$. Then, we define
		\begin{equation*}
			\begin{split}
				\hat{Y}_t := {\lambda}^{(\alpha),\varepsilon,\varrho}_t Y_t, \quad
				\hat{Z}_t := {\lambda}^{(\alpha),\varepsilon,\varrho}_t Z_t, \quad d\hat{M}_t := {\lambda}^{(\alpha),\varepsilon,\varrho}_t dM_t
			\end{split}
		\end{equation*}
		By applying the integration-by-parts formula, we obtain
		\begin{equation*}
			\begin{split}
				\hat{Y}_t = \hat{\xi} + \int_{t}^{T} \hat{f}\big(s,\hat{Y}_s,\hat{Z}_s\big) ds + \int_{t}^{T} \hat{g}\big(s,\hat{Y}_s\big) d\kappa_s - \int_{t}^T \hat{Z}_s d W_s- \int_{t}^T d\hat{M}_s,
			\end{split}
		\end{equation*}
		with
		\begin{equation*}
			\begin{split}
				\hat{\xi} &= {\lambda}^{(\alpha),\varepsilon,\varrho}_T \, \xi, \\
				\hat{f}(t,y,z) &= {\lambda}^{(\alpha),\varepsilon,\varrho}_t \, f\left(t, {\lambda}^{(-\alpha),-\varepsilon,-\varrho}_t y, {\lambda}^{(-\alpha),-\varepsilon,-\varrho}_t  z \right) - ({\alpha}_s+\varepsilon a^2_s) y,\\
				\hat{g}(t,y) &= {\lambda}^{(\alpha),\varepsilon,\varrho}_t \, g\left(t, {\lambda}^{(-\alpha),-\varepsilon,-\varrho}_t y \right) - \varrho y.
			\end{split}
		\end{equation*}
		Thus, if $(Y,Z)$ is a solution of the GBSDE (\ref{basic GBSDE}) associated with $(\xi, f,g,\kappa)$, then the process $(\hat{Y},\hat{Z})$ satisfies a similar GBSDE associated with $(\hat{\xi}, \hat{f}, \hat{g}, \kappa)$. The driver $\hat{f}$ satisfies the stochastic Lipschitz condition described in \textsc{(H3)} with the same stochastic process $(\eta_t)_{t \leq T}$. Moreover, if $f$ satisfies \textsc{(H2)}, it follows that the coefficient $\hat{f}$ also satisfies a similar monotonicity condition with the real-valued stochastic process $(\hat{\alpha}^\varepsilon_t)$ given by $\hat{\alpha}^\varepsilon_t = \alpha_t - (\alpha_t+\varepsilon a_t^2) = -\varepsilon a^2_t \leq 0$. Consequently, for any $t \in [0,T]$ and each $\varepsilon \geq 0$, we have $\hat{\alpha}^\varepsilon_t + \varepsilon a_t^2  = 0$. Therefore, this change of variable reduces the case to $\alpha_t + \varepsilon a_t^2 = 0$ for any $t \in [0,T]$ and for each given $\varepsilon \geq 0$. On the other hand, we can see that $\theta_t<0$ is not a severe restriction. Indeed, if $g$ satisfies \textsc{(H3)}, it follows that the coefficient $\hat{g}$ also satisfies a similar monotonicity condition with the real-valued stochastic process $(\widehat{\delta}_t)$ given by  $\widehat{\delta}_t = \theta_t - \varrho$. By choosing $\varrho$ large enough so that $\theta_t < \varrho$, we can always reduce the case where $\widehat{g}$ satisfies \textsc{(H3)} with $\theta$ negative.\\
		Finally, let $\varepsilon \geq 0$ and $(Y^\varepsilon, Z^\varepsilon, M^\varepsilon)$ be a solution of the GBSDE (\ref{basic GBSDE}) associated with $(\xi, f_\varepsilon, g, \kappa)$, where $f_\varepsilon$ satisfies \textsc{(H2)} with $(\alpha_t)$ replaced by $ (-\varepsilon a^2_t)$ and also verifies \textsc{(H4)}. Set $\lambda^{(\alpha),\varepsilon}_t=\exp\left\{-\int_{0}^{t}\alpha_s ds-\varepsilon A_t\right\}$ and define
		\begin{equation*}
			\begin{split}
				\hat{Y}^\varepsilon_t := \lambda^{(\alpha),\varepsilon}_t \, Y^\varepsilon_t, \quad
				\hat{Z}^\varepsilon_t := \lambda^{(\alpha),\varepsilon}_t \, Z^\varepsilon_t,\quad d\hat{M}^\varepsilon_t := \lambda^{(\alpha),\varepsilon}_t \, dM^\varepsilon_t.
			\end{split}
		\end{equation*}
		Using the integration-by-parts formula, we obtain
		\begin{equation*}
			\begin{split}
				\hat{Y}^\varepsilon_t = \hat{\xi}^\varepsilon + \int_{t}^{T} \hat{f}_\varepsilon\big(s,\hat{Y}^\varepsilon_s,\hat{Z}^\varepsilon_s\big) \, ds + \int_{t}^{T} \hat{g}_\varepsilon\big(s,\hat{Y}^\varepsilon_s\big) \, d\kappa_s - \int_{t}^T \hat{Z}^\varepsilon_s \, d W_s- \int_{t}^T d\hat{M}^\varepsilon_s
			\end{split}
		\end{equation*}
		with
		\begin{equation*}
			\begin{split}
				\hat{\xi}^\varepsilon &= \lambda^{(\alpha),\varepsilon}_T \, \xi, \\
				\hat{f}_\varepsilon(t,y,z) &= \lambda^{(\alpha),\varepsilon}_t \, f_\varepsilon\left(t, \lambda^{(-\alpha),-\varepsilon}_t y, \lambda^{(-\alpha),-\varepsilon}_t z \right) + (\alpha_t+\varepsilon a^2_t) \, y,\\
				\hat{g}_\varepsilon(t,y) &= \lambda^{(\alpha),\varepsilon}_t \, g\left(t, \lambda^{(-\alpha),-\varepsilon}_t y \right).
			\end{split}
		\end{equation*}
		Then, $\widehat{g}_\varepsilon$ satisfies \textsc{(H3)} with the negative process $(\theta_t)$, and $\widehat{f}_\varepsilon$ satisfies \textsc{(H2)} and \textsc{(H4)} with respect to the processes $(\alpha_t)$ and $(\eta_t)$, respectively.
	\end{remark}
	
	To facilitate the calculations, we will assume, for the remainder of this paper, that condition \textsc{(H2)} is satisfied with a process $(\alpha_t)_{t \leq T}$ such that $\alpha_t + \varepsilon a_t^2 = 0$ for each given $\varepsilon \geq 0$ and for all $t \in [0,T]$, and that \textsc{(H2)} holds for $\theta_t < 0$ for all $t \in [0,T]$. If this is not the case, the same change of variable described in Remark \ref{rmq essential} can be applied to reduce the situation to this case.
	
	We now provide the definition of an $\mathbb{L}^p$-solution for the GBSDE (\ref{basic GBSDE}):
	\begin{definition}\label{Def}
		Let $p \in (1, 2]$. A triplet $(Y_t, Z_t, M_t)_{t \leq T}$ is called an $\mathbb{L}^p$-solution of the GBSDE (\ref{basic GBSDE}) if the following conditions are satisfied:
		\begin{itemize}
			\item $(Y, Z, M)$ satisfies (\ref{basic GBSDE}).
			\item $(Y, Z, M)$ belongs to $\mathcal{E}^p_{\beta,\mu}$ for some $\beta, \mu > 0$.
		\end{itemize}
		
		Sometimes, we shall refer to the triplet $(Y_t, Z_t, M_t)_{t \leq T}$ as a solution in $\mathcal{E}^p_{\beta,\mu}$ if $(Y, Z, M)$ satisfies (\ref{basic GBSDE}) and belongs to $\mathcal{E}^p_{\beta,\mu}$ for some $\beta, \mu \geq 0$.
	\end{definition}

	Since the function $x \mapsto |x|^p$ is not sufficiently smooth for $p < 2$, It\^o's formula cannot be applied directly. Therefore, we need a generalization of Tanaka's formula for the multidimensional case. To this end, we introduce the notation $\check{x} = \frac{x}{|x|}\mathds{1}_{x \neq 0}$ for $x \in \mathbb{R}^d$. The following lemma extends the Meyer–It\^o formula as referenced in \cite{briand2003lp}. Although this result likely appears in earlier works, its proof is provided in \cite[Lemma 7]{kruse2016bsdes} (see also \cite[Lemma 2.2]{briand2003lp} for the Brownian case).
	
	In the context of the BSDE in a general filtration considered in \cite{kruse2016bsdes}, the additional generator term associated with $(\kappa_t)_{t \leq T}$ in our formulation (\ref{basic GBSDE}) is a continuous stochastic process, and our filtration supports a Brownian motion without requiring an independent jump measure. Consequently, the proof remains unchanged. The same observation applies if $\int_{0}^{\cdot} g(s,Y_s) \, d\kappa_s$ is replaced by an $\mathbb{R}$-valued, continuous, adapted process with locally integrable variation on $[0,T]$, denoted by $(G_t)_{t \leq T}$ with $G_0 = 0$. Therefore, we omit the proof for brevity.
	\begin{lemma}\label{Lp basic decomposition}
		Let $(F_t)_{t \leq T}$ and $(Z_t)_{t \leq T}$ be two progressively measurable processes with values respectively in $\mathbb{R}^d$ and $\mathbb{R}^{d \times k}$ such that $\mathbb{P}$-a.s.
		$$
		\int_{0}^{T}\left\{\left(\left|F_s\right|+\left\|Z_s\right\|^2\right)ds+d\left[M\right]_s\right\}<+\infty.
		$$
		We consider the $\mathbb{R}^d$-valued semimartingale $(X_t)_{t \leq T}$ defined by 
		\begin{equation}\label{semimartingale decom}
			X_t=X_0+\int_{0}^{t}F_s ds+\int_0^t dG_s+\int_{0}^{t}Z_s dW_s+\int_{0}^{t}dM_s
		\end{equation}
		Then, for any $p \geq 1$, the process $(\big|X_t\big|^p)_{t \leq T}$ is an $\mathbb{R}$-semimartingale with the following decomposition 
		\begin{equation*}
			\begin{split}
				&\left|X_t\right|^p\\
				&=\left|X_0\right|^p+\frac{1}{2}\mathds{1}_{p=1}L_t+p\int_{0}^{t} \left|X_s\right|^{p-1} \check{X}_s F_s ds+p\int_{0}^{t} \left|X_s\right|^{p-1} \check{X}_s dG_s \\
				&+p\int_{0}^{t} \left|X_s\right|^{p-1} \check{X}_s Z_s dW_s+p\int_{0}^{t} \left|X_{s-}\right|^{p-1} \check{X}_{s-} dM_s\\
				&+\sum_{0 <s \leq t} \left\{\left|X_{s-}+\Delta M_s\right|^p-\left|X_{s-}\right|^p-p\left|X_{s-}\right|^{p-1}\check{X}_{s-}\Delta M_s\right\}\\
				&+\frac{p}{2}\int_{0}^{t}\left|X_s\right|^{p-2}\mathds{1}_{X_s \neq 0}\left\{(2-p)\left(\left\|Z_s\right\|^2 -(\check{X}_s)^\ast Z_s Z^\ast_s \check{X}_s\right)  +(p-1)\left\|Z_s\right\|^2\right\}ds 
				\\
				&+\frac{p}{2}\int_{0}^{t}\left|X_s\right|^{p-2}\mathds{1}_{X_s \neq 0}\left\{(2-p)\left(d\left[M\right]^c_s-(\check{X}_s)^\ast d\left[M,M\right]^c_s \check{X}_s\right)+(p-1)d\left[M\right]^c_s \right\},
			\end{split}
		\end{equation*}
		where $(L_t)_{t \leq T}$ is a continuous, increasing process with $L_0=0$, which increases only on the boundary of the random set $\{t \in [0,T] : X_{t-}=X_t=0\}$.
	\end{lemma}

	Let $(X_t)_{t \leq T}$ be an $\mathbb{R}^d$-valued semimartingale with the decomposition (\ref{semimartingale decom}), and let $(\bar{\kappa}_t)_{t \leq T}$ be a continuous adapted process with locally integrable variation on $[0,T]$. Let $A$ be the process defined in assumption \textbf{(H-M)$_\mathbf{p}$}. Using Lemma \ref{Lp basic decomposition} and the integration-by-parts formula (see Corollary 2 in \cite[p. 68]{protter2005stochastic}) for the product of semimartingales $\big(e^{\frac{p}{2}\beta A_t +\frac{p}{2} \mu \bar{\kappa}_t} |X_t|^p\big)_{t \leq T}$, we obtain the following result:
	\begin{corollary}\label{basic coro}
		For any $p \in (1,2)$ and all $\beta, \mu >0$, we have
		\begin{equation*}
			\begin{split}
				&e^{\frac{p}{2}\beta A_t +\frac{p}{2} \mu \bar{\kappa}_t}\left|X_t\right|^p\\
				=&\left|X_0\right|^p+\frac{p}{2}\beta \int_{0}^{t}e^{\frac{p}{2}\beta A_s +\frac{p}{2} \mu \bar{\kappa}_s}\left|X_s\right|^pdA_s+\frac{p}{2}\mu  \int_{0}^{t}e^{\frac{p}{2}\beta A_s +\frac{p}{2} \mu \bar{\kappa}_s}\left|X_s\right|^pd\bar{\kappa}_s\\
				&+p\int_{0}^{t}e^{\frac{p}{2}\beta A_s +\frac{p}{2} \mu \bar{\kappa}_s} \left|X_s\right|^{p-1} \check{X}_s F_s ds+p\int_{0}^{t} e^{\frac{p}{2}\beta A_s +\frac{p}{2} \mu \bar{\kappa}_s}\left|X_s\right|^{p-1} \check{X}_s d G_s\\
				&+p\int_{0}^{t}e^{\frac{p}{2}\beta A_s +\frac{p}{2} \mu \bar{\kappa}_s} \left|X_s\right|^{p-1} \check{X}_s Z_s dW_s+p\int_{0}^{t}e^{\frac{p}{2}\beta A_s +\frac{p}{2} \mu \bar{\kappa}_s} \left|X_{s-}\right|^{p-1} \check{X}_{s-} dM_s\\
				&+\sum_{0 <s \leq t} e^{\frac{p}{2}\beta A_s +\frac{p}{2} \mu \bar{\kappa}_s}\left\{\left|X_{s-}+\Delta M_s\right|^p-\left|X_{s-}\right|^p-p\left|X_{s-}\right|^{p-1}\check{X}_{s-}\Delta M_s\right\}\\
				&+\frac{p}{2}\int_{0}^{t}e^{\frac{p}{2}\beta A_s +\frac{p}{2} \mu \bar{\kappa}_s}\left|X_s\right|^{p-2}\mathds{1}_{X_s \neq 0} \left\{(2-p) \left(\left\|Z_s\right\|^2-(\check{X}_s)^\ast Z_s Z^\ast_s \check{X}_s\right)\right.\\
				&\left.+(p-1)\left\|Z_s\right\|^2\right\}ds+\frac{p}{2}\int_{0}^{t}e^{\frac{p}{2}\beta A_s +\frac{p}{2} \mu \bar{\kappa}_s}\left|X_s\right|^{p-2}\mathds{1}_{X_s \neq 0}\left\{(2-p)\left(d\left[M\right]^c_s \right.\right.\\
				&\left.\left.-(\check{X}_s)^\ast d\left[M,M\right]^c_s \check{X}_s\right)+(p-1)d\left[M\right]^c_s \right\}.
			\end{split}
		\end{equation*}
	\end{corollary}

	Let $p \in (1,2)$, and let $(Y^1_t, Z^1_t)_{t \leq T}$ and $(Y^2_t, Z^2_t)_{t \leq T}$ be two $\mathbb{L}^p$-solutions of the GBSDE (\ref{basic GBSDE}) associated with the data $(\xi^1, f^1, g^1, \kappa^1)$ and $(\xi^2, f^2, g^2, \kappa^2)$, respectively. Define $\widehat{\mathcal{R}} = \mathcal{R}^1 - \mathcal{R}^2$, where $\mathcal{R} \in \left\{Y, Z, M, \kappa \right\}$. Then, for all $0 \leq t \leq t^{\prime} \leq T$, the pair $(\widehat{Y}_t, \widehat{Z}_t, \widehat{M}_t)_{t \leq T}$ satisfies the following GBSDE:
	\begin{equation}\label{Khdj}
		\widehat{Y}_t = \widehat{Y}_{t^{\prime}} + \int_{t}^{t^{\prime}} F_s \, ds + \int_{t}^{t^{\prime}} dG_s - \int_{t}^{t^{\prime}} \widehat{Z}_s \, dW_s- \int_{t}^{t^{\prime}} d\widehat{M}_s ,
	\end{equation}
	where $F_s = f^1(s, Y^1_s, Z^1_s) - f^1(s, Y^2_s, Z^2_s)$ and $dG_s = g^1(s, Y^1_s) \, d\kappa^1_s - g^2(s, Y^2_s) \, d\kappa^2_s$. Note that the process $(G_t)_{t \leq T}$ can be written, in the sense of signed measures on $[0,T]$, as follows:
	$$
	dG_s = \big(g^1(s, Y^1_s) - g^1(s, Y^2_s) \big) \, d\kappa^1_s + g^1(s, Y^2_s) \, d\widehat{\kappa}_s + \big(g^1(s, Y^2_s) - g^2(s, Y^2_s) \big) \, d\kappa^2_s.
	$$
	Next, let $\bar{\kappa}_t := \|\widehat{\kappa}\|_t + \kappa^2_t$ for $t \in [0, T]$ and $\beta, \mu > 0$. Then, we derive the following:
	\begin{corollary}\label{used in a priori estimates}
		Let $(Y^1_t, Z^1_t, M^1_t)_{t \leq T}$ and $(Y^2_t, Z^2_t, M^2_t)_{t \leq T}$ be two $\mathbb{L}^p$-solutions of the GBSDE (\ref{basic GBSDE}) associated with the data $(\xi^1, f^1, g^1, \kappa^1)$ and $(\xi^2, f^2, g^2, \kappa^2)$, respectively. Set $c(p) = \frac{p(p-1)}{2}$ and let $\tau \in \mathcal{T}_{0,T}$, then we have:
		\begin{equation}
				\begin{split}
					& e^{\frac{p}{2}\beta A_{t \wedge \tau} +\frac{p}{2} \mu \bar{\kappa}_{t \wedge \tau}}\big|\widehat{Y}_{t \wedge \tau}\big|^p+\frac{p}{2}\beta \int_{t \wedge \tau}^{\tau}e^{\frac{p}{2}\beta A_{s} +\frac{p}{2} \mu \bar{\kappa}_{s}}\big|\widehat{Y}_s\big|^pdA_s\\
					&+\frac{p}{2}\mu  \int_{t \wedge \tau}^{\tau}e^{\frac{p}{2}\beta A_{s} +\frac{p}{2} \mu \bar{\kappa}_{s}}\big|\widehat{Y}_s\big|^p(d\|\widehat{\kappa}\|_s+d\kappa^2_s)\\
					&+c(p)\int_{t \wedge \tau}^{\tau}e^{\frac{p}{2}\beta A_{s} +\frac{p}{2} \mu \bar{\kappa}_{s}}\big|\widehat{Y}_s\big|^{p-2}\mathds{1}_{Y_s \neq 0}\big\|\widehat{Z}_s\big\|^2ds\\
					&+c(p)\int_{t \wedge \tau}^{\tau}e^{\frac{p}{2}\beta A_{s} +\frac{p}{2} \mu \bar{\kappa}_{s}}\big|\widehat{Y}_s\big|^{p-2}\mathds{1}_{Y_s \neq 0}d \big[\widehat{M}\big]^c_s\\
					\leq &\bar{\Theta}^{\frac{p}{2}\beta ,\frac{p}{2}\mu}_\tau\big|\widehat{Y}_\tau\big|^p
					+p\int_{t \wedge \tau}^{\tau}e^{\frac{p}{2}\beta A_{s} +\frac{p}{2} \mu \bar{\kappa}_{s}} \big|\widehat{Y}_s\big|^{p-1} \check{\widehat{Y}}_s \big(f^1(s,Y^1_s,Z^1_s)-f^1(s,Y^2_s,Z^2_s)\big) ds\\
					&+p\int_{t \wedge \tau}^{\tau} e^{\frac{p}{2}\beta A_{s} +\frac{p}{2} \mu \bar{\kappa}_{s}}\big|\widehat{Y}_s\big|^{p-1} \check{\widehat{Y}}_s \left\{\big(g^1(s,Y^1_s)-g^1(s,Y^2_s)\big)d\kappa^1_s+g^1(s,Y^2_s)d\widehat{\kappa}_s\right\}\\
					&+p\int_{t \wedge \tau}^{\tau} e^{\frac{p}{2}\beta A_{s} +\frac{p}{2} \mu \bar{\kappa}_{s}}\big|\widehat{Y}_s\big|^{p-1} \check{\widehat{Y}}_s\big(g^1(s,Y^2_s)-g^2(s,Y^2_s)\big)d\kappa^2_s\\
					&-p\int_{t \wedge \tau}^{\tau}e^{\frac{p}{2}\beta A_{s} +\frac{p}{2} \mu \bar{\kappa}_{s}} \big|\widehat{Y}_s\big|^{p-1} \check{\widehat{Y}}_s \widehat{Z}_s dW_s-p\int_{t \wedge \tau}^{\tau}\bar{\Theta}^{\frac{p}{2}\beta ,\frac{p}{2}\mu}_s \big|\widehat{Y}_{s-}\big|^{p-1} \check{\widehat{Y}}_{s-} d\widehat{M}_s\\
					&-\sum_{t \wedge \tau <s \leq \tau} e^{\frac{p}{2}\beta A_{s} +\frac{p}{2} \mu \bar{\kappa}_{s}}\left\{\big|\widehat{Y}_{s-}+\Delta \widehat{M}_s\big|^p-\big|\widehat{Y}_{s-}\big|^p-p\big|\widehat{Y}_{s-}\big|^{p-1}\check{\widehat{Y}}_{s-}\Delta \widehat{M}_s\right\}.
				\end{split}
		\end{equation}
	\end{corollary}

	\begin{proof}
		First, using the Cauchy-Schwarz inequality, we have $\mathbb{P}$-a.s.
		\begin{equation}\label{for Z}
			(2 - p) \big(\big\| \widehat{Z}_s \big\|^2 - (\check{X}_s)^\ast \widehat{Z}_s \widehat{Z}^\ast_s \check{X}_s\big) \leq 0
		\end{equation}
		Next, recall that $[M, M] = \left\{\left[M^i, M^j\right]^c; 1 \leq i, j \leq d\right\}$. We define $\mathcal{X}^{\beta, \mu}_s = e^{\frac{p}{2}\beta A_{s} +\frac{p}{2} \mu \bar{\kappa}_{s}} \left| X_s \right|^{p-4} \mathds{1}_{X_s \neq 0}$. By applying the Kunita-Watanabe inequality (as in Remark \ref{Good rmq too}), we have $\mathbb{P}$-a.s.
		\begin{equation}\label{for M}
			\begin{split}
				&\int_{t \wedge \tau}^{\tau}e^{\frac{p}{2}\beta A_{s} +\frac{p}{2} \mu \bar{\kappa}_{s}}\left|X_s\right|^{p-2}\mathds{1}_{X_s \neq 0}(\check{X}_s)^\ast d\left[M,M\right]^c_s \check{X}_s\\
				&=\sum_{i,j=1}^d\int_{t \wedge \tau}^{\tau}\bar{\Theta}^{\frac{p}{2}\beta ,\frac{p}{2}\mu}_s\left|X_s\right|^{p-4}\mathds{1}_{X_s \neq 0}X^i_s X^j_s d\left[M^i,M^j\right]^c_s \\
				& \leq \sum_{i,j=1}^d\left(\int_{t \wedge \tau}^{\tau}\mathcal{X}^{\beta,\mu}_s\left|X^i_s\right|^2d\left[M^i\right]^c_s \right)^{\frac{1}{2}}\left(\int_{0}^{t}\mathcal{X}^{\beta,\mu}_s \left|X^j_s\right|^2 d\left[M^j\right]^c_s \right)^{\frac{1}{2}}\\
				& \leq \frac{1}{2} \left(\sum_{i=1}^d\int_{t \wedge \tau}^{\tau}\mathcal{X}^{\beta,\mu}_s\left|X^i_s\right|^2d\left[M^i\right]^c_s +\sum_{j=1}^d\int_{t \wedge \tau}^{\tau}\mathcal{X}^{\beta,\mu}_s \left|X^j_s\right|^2 d\left[M^j\right]^c_s \right) \\
				& \leq \int_{t \wedge \tau}^{\tau}e^{\frac{p}{2}\beta A_{s} +\frac{p}{2} \mu \bar{\kappa}_{s}}\left|X_s\right|^{p-2}\mathds{1}_{X_s \neq 0}d\left[M\right]^c_s 
			\end{split}
		\end{equation}
		Finally, using Corollary \ref{basic coro} with $X = \widehat{Y}$, and applying (\ref{Khdj}) along with (\ref{for Z})-(\ref{for M}), completes the proof.
	\end{proof}

	From this point forward and throughout the remainder of the paper, unless otherwise specified, we assume $p \in (1, 2)$.
	
	\section{A priori estimates and uniqueness}
	\label{sec3}
	Let $(Y^1_t, Z^1_t, M^1_t)_{t \leq T}$ and $(Y^2_t, Z^2_t, M^2_t)_{t \leq T}$ be two $\mathbb{L}^p$-solutions of the GBSDE (\ref{basic GBSDE}) associated with the data $(\xi^1, f^1, g^1, \kappa^1)$ and $(\xi^2, f^2, g^2, \kappa^2)$, respectively, satisfying condition \textbf{(H-M)$_\mathbf{p}$}. Define $\mathcal{R} = \mathcal{R}^1 - \mathcal{R}^2$ with $\mathcal{R} \in \left\{Y, Z, M, \xi, f, g,\kappa \right\}$. Then, we have
	\begin{proposition}\label{propo1}
		For any $\beta, \mu > \frac{2(p-1)}{p}$, there exists a constant  $\mathfrak{c}_{\beta,\mu,p,\epsilon}$ such that
		\begin{equation*}
			\begin{split}
				&\mathbb{E}\left[\sup_{t \in [0,T]}e^{\frac{p}{2}\beta A_{t} +\frac{p}{2} \mu \bar{\kappa}_{t}}\big|\widehat{Y}_t\big|^p\right]+\mathbb{E} \int_{0}^{T}e^{\frac{p}{2}\beta A_{s} +\frac{p}{2} \mu \bar{\kappa}_{s}}\big|\widehat{Y}_s\big|^p (dA_s+d\|\widehat{\kappa}\|_s+d\kappa^2_s)\\
				&+\mathbb{E}\left[ \left( \int_{0}^{T}e^{\beta A_{s} +\mu \bar{\kappa}_{s}}\big\|\widehat{Z}_s\big\|^2ds\right)^{\frac{p}{2}}\right]  
				+\mathbb{E}\left[ \left( \int_{0}^{T}e^{\beta A_{s} + \mu \bar{\kappa}_{s}}d\big[\widehat{M}\big]_s\right)^{\frac{p}{2}}\right] \\
				\leq & \mathfrak{c}_{\beta,\mu,p,\epsilon}\left( \mathbb{E}\left[ e^{\frac{p}{2}\beta A_{T} +\frac{p}{2} \mu \bar{\kappa}_{T}}\big|\widehat{\xi}\big|^p\right] 
				+\int_{0}^{T}e^{\beta A_{s} + \mu \bar{\kappa}_{s}} \big|\widehat{f}(s,Y^2_s,Z^2_s)\big|^pds\right. \\
				&\left. +\mathbb{E}\int_{0}^{T} e^{\beta A_{s} + \mu \bar{\kappa}_{s}} \big|g^1(s,Y^2_s)\big|^pd\|\widehat{\kappa}\|_s +\mathbb{E}\int_{0}^{T} e^{\beta A_{s} + \mu \bar{\kappa}_{s}}\big|\widehat{g}(s,Y^2_s)\big|^pd\kappa^2_s\right).
			\end{split}
		\end{equation*}
	\end{proposition}

	\begin{proof}
		As in \cite[Proposition 1]{elmansourielotmani2024}, it suffices to prove the result in the case where $\|\widehat{\kappa}\|_T + \kappa^2_T$ is a bounded random variable, and then apply Fatou's Lemma.\\
		Using assumptions (H2)–(H4) on the drivers $f$ and $g$, along with the basic inequality $ab \leq \frac{a^2}{2\varepsilon} + \frac{\varepsilon b^2}{2}$ for all $\varepsilon > 0$ and Remark \ref{rmq essential}, we derive the following:
		\begin{equation}\label{CV}
			\begin{split}
				&\widehat{Y}_s\left(f^1(s,Y^1_s,Z^1_s)-f^2(s,Y^2_s,Z^2_s)\right)
				\leq  \frac{c(p)}{2} \big\|\widehat{Z}_s\big\|^2+\big|\widehat{Y}_s\big| \big|\widehat{f}(s,Y^2_s,Z^2_s)\big|
			\end{split}
		\end{equation}
		and
		\begin{equation}\label{CV1}
			\begin{split}
				&\widehat{Y}_s\left(g^1(s,Y^1_s)-g^1(s,Y^2_s)\right)d\kappa^1_s+\widehat{Y}_sg^1(s,Y^2_s)d\widehat{\kappa}_s
				+\widehat{Y}_s\big(g^1(s,Y^2_s)-g^2(s,Y^2_s)\big)d\kappa^2_s\\
				& \leq\big|\widehat{Y}_s\big| \big|g^1(s,Y^2_s)\big|d\|\widehat{\kappa}\|_s+\big|\widehat{Y}_s\big| \big|\widehat{g}(s,Y^2_s)\big|d\kappa^2_s
			\end{split}
		\end{equation}
		Next, using (\ref{CV})–(\ref{CV1}) along with Corollary \ref{used in a priori estimates} on $[t, T]$ with $\tau = T$ and the fact that $\big|\widehat{Y}_s\big|^{p-1} \check{\widehat{Y}}_s = \big|\widehat{Y}_s\big|^{p-2} \mathds{1}_{Y_s \neq 0} \, \widehat{Y}_s$, we obtain
		\begin{equation}\label{Coming back to this Itos formula}
			\begin{split}
				&e^{\frac{p}{2}\beta A_{t} +\frac{p}{2} \mu \bar{\kappa}_{t}}\big|\widehat{Y}_{t}\big|^p+\frac{p}{2}\beta \int_{t }^{T}e^{\frac{p}{2}\beta A_{s} +\frac{p}{2} \mu \bar{\kappa}_{s}}\big|\widehat{Y}_s\big|^pdA_s\\
				&+\frac{p}{2}\mu  \int_{t  }^{T}e^{\frac{p}{2}\beta A_{s} +\frac{p}{2} \mu \bar{\kappa}_{s}}\big|\widehat{Y}_s\big|^p(d\|\widehat{\kappa}\|_s+d\kappa^2_s)\\
				&+\frac{c(p)}{2}\int_{t}^{T}e^{\frac{p}{2}\beta A_{s} +\frac{p}{2} \mu \bar{\kappa}_{s}}\big|\widehat{Y}_s\big|^{p-2}\mathds{1}_{Y_s \neq 0}\big\|\widehat{Z}_s\big\|^2ds\\
				&+c(p)\int_{t }^{T}e^{\frac{p}{2}\beta A_{s} +\frac{p}{2} \mu \bar{\kappa}_{s}}\big|\widehat{Y}_s\big|^{p-2}\mathds{1}_{Y_s \neq 0}d \big[\widehat{M}\big]^c_s\\
				\leq &e^{\frac{p}{2}\beta A_{T} +\frac{p}{2} \mu \bar{\kappa}_{T}}\big|\widehat{\xi}\big|^p
				+p\int_{t}^{T}e^{\frac{p}{2}\beta A_{s} +\frac{p}{2} \mu \bar{\kappa}_{s}} \big|\widehat{Y}_s\big|^{p-1}  \big|\widehat{f}(s,Y^2_s,Z^2_s)\big|ds\\
				&+p\int_{t}^{T} e^{\frac{p}{2}\beta A_{s} +\frac{p}{2} \mu \bar{\kappa}_{s}}\big|\widehat{Y}_s\big|^{p-1} \left\{\big|g^1(s,Y^2_s)\big|d\|\widehat{\kappa}\|_s+\big|\widehat{g}(s,Y^2_s)\big|d\kappa^2_s\right\}\\
				&-p\int_{t}^{T}e^{\frac{p}{2}\beta A_{s} +\frac{p}{2} \mu \bar{\kappa}_{s}} \big|\widehat{Y}_s\big|^{p-1} \check{\widehat{Y}}_s \widehat{Z}_s dW_s-p\int_{t}^{T}e^{\frac{p}{2}\beta A_{s} +\frac{p}{2} \mu \bar{\kappa}_{s}} \big|\widehat{Y}_{s-}\big|^{p-1} \check{\widehat{Y}}_{s-} d\widehat{M}_s\\
				&-\sum_{t <s \leq T} e^{\frac{p}{2}\beta A_{s} +\frac{p}{2} \mu \bar{\kappa}_{s}}\left\{\big|\widehat{Y}_{s-}+\Delta \widehat{M}_s\big|^p-\big|\widehat{Y}_{s-}\big|^p-p\big|\widehat{Y}_{s-}\big|^{p-1}\check{\widehat{Y}}_{s-}\Delta \widehat{M}_s\right\}.
			\end{split}
		\end{equation}
		By applying H\"older's inequality $\int \big|h\big|^{p-1} \big|\ell\big| \, d \vartheta \leq \left( \int \big|h\big|^p \, d \vartheta \right)^{\frac{p-1}{p}} \left( \int \big|\ell\big|^p \, d \vartheta \right)^{\frac{1}{p}}$, Young's inequality $a^{\frac{p-1}{p}} b^{\frac{1}{p}} \leq \frac{p-1}{p} \, a + \frac{1}{p} \, b$, and assumption (H6), we derive
		\begin{equation}\label{fes}
				\begin{split}
					&p\int_{t}^{T}e^{\frac{p}{2}\beta A_{s} +\frac{p}{2} \mu \bar{\kappa}_{s}} \big|\widehat{Y}_s\big|^{p-1}  \big|\widehat{f}(s,Y^2_s,Z^2_s)\big|ds\\
					=&p\int_{t}^{T} \left(e^{\frac{p-1}{2}\beta A_{s} +\frac{p-1}{2} \mu \bar{\kappa}_{s}}\big|\widehat{Y}_s\big|^{p-1} a^{\frac{2(p-1)}{p}}_s\right) \left(e^{\frac{\beta}{2}A_{s} +\frac{\mu}{2}  \bar{\kappa}_{s}}\big|\widehat{f}(s,Y^2_s,Z^2_s)\big|a^{\frac{2(1-p)}{p}}_s\right)ds\\
					\leq& (p-1) \int_{t}^{T} e^{\frac{p}{2}\beta A_{s} +\frac{p}{2} \mu \bar{\kappa}_{s}}\big|\widehat{Y}_s\big|^{p} dA_s +\int_{t}^{T} e^{\frac{p}{2}\beta A_{s} +\frac{p}{2} \mu \bar{\kappa}_{s}}\big|\widehat{f}(s,Y^2_s,Z^2_s)\big|^p a^{2(1-p)}_s ds\\
					\leq& (p-1) \int_{t}^{T} e^{\frac{p}{2}\beta A_{s} +\frac{p}{2} \mu \bar{\kappa}_{s}}\big|\widehat{Y}_s\big|^{p} dA_s +\frac{1}{\epsilon^{2(p-1)}} \int_{t}^{T} e^{\beta A_{s} + \mu \bar{\kappa}_{s}}\big|\widehat{f}(s,Y^2_s,Z^2_s)\big|^p ds
				\end{split}
		\end{equation}
		Using a similar argument, we obtain
		\begin{equation}\label{ges}
			\begin{split}
				&p\int_{t}^{T} e^{\frac{p}{2}\beta A_{s} +\frac{p}{2} \mu \bar{\kappa}_{s}}\big|\widehat{Y}_s\big|^{p-1} \left\{\big|g^1(s,Y^2_s)\big|d\|\widehat{\kappa}\|_s+\big|g^1(s,Y^2_s)-g^2(s,Y^2_s)\big|d\kappa^2_s\right\}\\
				&\leq (p-1)\int_{t}^{T} e^{\frac{p}{2}\beta A_{s} +\frac{p}{2} \mu \bar{\kappa}_{s}}\big|\widehat{Y}_s\big|^{p}\left\{d\|\widehat{\kappa}\|_s+d\kappa^2_s\right\}\\
				&\qquad+ \int_{t}^{T} e^{\beta A_{s} + \mu \bar{\kappa}_{s}}\big|g^1(s,Y^2_s)\big|^p d\|\widehat{\kappa}\|_s+ \int_{t}^{T} e^{\beta A_{s} + \mu \bar{\kappa}_{s}}\big|\widehat{g}(s,Y^2_s)\big|^p d\kappa^2_s
			\end{split}
		\end{equation}
		Next, from \cite[Lemma 9]{kruse2016bsdes}, the jump part of the quadratic variation is controlled by a non-decreasing process involving the jumps of $Y$, as described below
		\begin{equation}\label{jumps for M}
				\begin{split}
					&\sum_{t <s \leq T} e^{\frac{p}{2}\beta A_{s} +\frac{p}{2} \mu \bar{\kappa}_{s}}
					\Bigg\{
					\big|\widehat{Y}_{s-}+\Delta \widehat{M}_s\big|^p
					-\big|\widehat{Y}_{s-}\big|^p
					-p\big|\widehat{Y}_{s-}\big|^{p-1}\check{\widehat{Y}}_{s-}\Delta \widehat{M}_s
					\Bigg\}\\
					&\geq c(p)\sum_{t <s \leq T} e^{\frac{p}{2}\beta A_{s} +\frac{p}{2} \mu \bar{\kappa}_{s}}
					\big|\Delta \widehat{M}_s\big|^2 
					\left(\big|\widehat{Y}_{s-}\big|^2 \vee \big|\widehat{Y}_{s-}+\Delta \widehat{M}_s\big|^2\right)^{\frac{p-2}{2}}
					\mathds{1}_{\left|\widehat{Y}_{s-}\right| \vee \left|\widehat{Y}_{s-}+\Delta \widehat{M}_s\right| \neq 0}.
				\end{split}
		\end{equation}
		Additionally, from the dynamics of the process $\widehat{Y}$ given by (\ref{Khdj}), we know that $\Delta \widehat{Y}_s = \Delta \widehat{M}_s$. Therefore, $\widehat{Y}_s = \widehat{Y}_{s-} + \Delta \widehat{M}_s$ for any $s \in [0, T]$.\\
		Returning to (\ref{Coming back to this Itos formula}) and using (\ref{fes}), (\ref{ges}), and (\ref{jumps for M}), we have
		\begin{equation}\label{back2}
			\begin{split}
				&e^{\frac{p}{2}\beta A_{t} +\frac{p}{2} \mu \bar{\kappa}_{t}}\big|\widehat{Y}_{t}\big|^p+ \left(\frac{p}{2}\beta-(p-1)\right)\int_{t}^{T}e^{\frac{p}{2}\beta A_{s} +\frac{p}{2} \mu \bar{\kappa}_{s}}\big|\widehat{Y}_s\big|^pdA_s\\
				+&\left(\frac{p}{2}\beta-(p-1)\right) \int_{t  }^{T}e^{\frac{p}{2}\beta A_{s} +\frac{p}{2} \mu \bar{\kappa}_{s}}\big|\widehat{Y}_s\big|^p(d\|\widehat{\kappa}\|_s+d\kappa^2_s)\\
				+&\frac{c(p)}{2}\int_{t}^{T}e^{\frac{p}{2}\beta A_{s} +\frac{p}{2} \mu \bar{\kappa}_{s}}\big|\widehat{Y}_s\big|^{p-2}\mathds{1}_{Y_s \neq 0}\big\|\widehat{Z}_s\big\|^2ds\\
				+&c(p)\int_{t}^{T}e^{\frac{p}{2}\beta A_{s} +\frac{p}{2} \mu \bar{\kappa}_{s}}\big|\widehat{Y}_s\big|^{p-2}\mathds{1}_{Y_s \neq 0}d \big[\widehat{M}\big]^c_s\\
				+&c(p)\sum_{t <s \leq T} e^{\frac{p}{2}\beta A_{s} +\frac{p}{2} \mu \bar{\kappa}_{s}} \big|\Delta \widehat{M}_s\big|^2 \left(\big|\widehat{Y}_{s-}\big|^2 \vee \big|\widehat{Y}_{s}\big|^2\right)^{\frac{p-2}{2}} \mathds{1}_{\left|\widehat{Y}_{s-}\right| \vee \left|\widehat{Y}_{s}\right| \neq 0}\\
				\leq &  e^{\frac{p}{2}\beta A_{T} +\frac{p}{2} \mu \bar{\kappa}_{T}}\big|\widehat{\xi}\big|^p
				+\int_{t}^{T}e^{\beta A_{s} + \mu \bar{\kappa}_{s}} \big|\widehat{f}(s,Y^2_s,Z^2_s)\big|^pds \\
				& +\int_{t}^{T} e^{\beta A_{s} + \mu \bar{\kappa}_{s}} \big|g^1(s,Y^2_s)\big|d\|\widehat{\kappa}\|_s +\int_{t}^{T} e^{\beta A_{s} + \mu \bar{\kappa}_{s}}\big|\widehat{g}(s,Y^2_s)\big|^pd\kappa^2_s\\
				&-p\int_{t}^{T}e^{\frac{p}{2}\beta A_{s} +\frac{p}{2} \mu \bar{\kappa}_{s}} \big|\widehat{Y}_s\big|^{p-1} \check{\widehat{Y}}_s \widehat{Z}_s dW_s-p\int_{t}^{T}e^{\frac{p}{2}\beta A_{s} +\frac{p}{2} \mu \bar{\kappa}_{s}} \big|\widehat{Y}_{s-}\big|^{p-1} \check{\widehat{Y}}_{s-} d\widehat{M}_s.
			\end{split}
		\end{equation}
		Let us set 
		$$\Lambda:=\int_{0}^{\cdot}e^{\frac{p}{2}\beta A_{s} +\frac{p}{2} \mu \bar{\kappa}_{s}} \big|\widehat{Y}_s\big|^{p-1} \check{\widehat{Y}}_s \widehat{Z}_s dW_s$$ 
		\mbox{and } 
		$$\Xi:=\int_{0}^{\cdot}e^{\frac{p}{2}\beta A_{s} +\frac{p}{2} \mu \bar{\kappa}_{s}} \big|\widehat{Y}_{s-}\big|^{p-1} \check{\widehat{Y}}_{s-} d\widehat{M}_s.$$ 
		It follows from the Burkholder-Davis-Gundy inequality (BDG; see, e.g, Theorem 48 in \cite[p. 193]{protter2005stochastic}) that the local martingales $\Lambda$ and $\Xi$ are uniformly integrable martingales with zero expectation.\\
		Indeed, for the continuous local martingale part $\Lambda$, we have, by Young's inequality,
		\begin{equation}\label{L2}
			\begin{split}
				\mathbb{E}\left[\sup_{t \in [0,T]} \big|\Lambda_t\big|\right]&
				\leq \mathfrak{c} \mathbb{E}\left[\left(\int_{0}^{T}e^{p\beta A_{s} +p \mu \bar{\kappa}_{s}} \big|\widehat{Y}_s\big|^{2(p-1)} \|\widehat{Z}_s\|^2 ds\right)^{\frac{1}{2}}\right]\\
				& \leq \mathfrak{c} \mathbb{E}\left[\left(\sup_{t \in [0,T]}e^{\frac{p-1}{2}\beta A_{t} +\frac{p-1}{2} \mu \bar{\kappa}_{t}} \big|\widehat{Y}_t\big|^{p-1}\right) \left(\int_{0}^{T}e^{\beta A_{s} + \mu \bar{\kappa}_{s}}  \|\widehat{Z}_s\|^2 ds\right)^{\frac{1}{2}}\right]\\
				&\leq \frac{(p-1)\mathfrak{c}}{p}\mathbb{E}\left[\sup_{t \in [0,T]}e^{\frac{p}{2}\beta A_{t} +\frac{p}{2} \mu \bar{\kappa}_{t}} \big|\widehat{Y}_t\big|^{p}\right]+\frac{\mathfrak{c}}{p}\mathbb{E}\left[\left(\int_{0}^{T}e^{\beta A_{s} + \mu \bar{\kappa}_{s}}  \|\widehat{Z}_s\|^2 ds\right)^{\frac{p}{2}}\right]\\
				&<+\infty.
			\end{split}
		\end{equation}
		A similar argument holds for the RCLL local martingale part $\Xi$, where we have
		\begin{equation}\label{L3}
			\begin{split}
				\mathbb{E}\left[\sup_{t \in [0,T]} \big|\Xi_t\big|\right]
				& \leq \mathfrak{c} \mathbb{E}\left[\left(\int_{0}^{T}e^{p\beta A_{s} +p \mu \bar{\kappa}_{s}} \big|\widehat{Y}_{s-}\big|^{2(p-1)} d \big[\widehat{ M}\big]_s\right)^{\frac{1}{2}}\right]\\
				& \leq \mathfrak{c} \mathbb{E}\left[\left(\sup_{t \in [0,T]}e^{\frac{p-1}{2}\beta A_{s} +\frac{p-1}{2} \mu \bar{\kappa}_{s}} \big|\widehat{Y}_t\big|^{p-1}\right)
				  \left(\int_{0}^{T}e^{\beta A_{s} + \mu \bar{\kappa}_{s}}  d \big[\widehat{ M}\big]_s\right)^{\frac{1}{2}}\right]\\
				&\leq \frac{(p-1)\mathfrak{c}}{p}\mathbb{E}\left[\sup_{t \in [0,T]}e^{\frac{p}{2}\beta A_{t} +\frac{p}{2} \mu \bar{\kappa}_{t}} \big|\widehat{Y}_t\big|^{p}\right]+\frac{\mathfrak{c}}{p}\mathbb{E}\left[\left(\int_{0}^{T}e^{\beta A_{s} + \mu \bar{\kappa}_{s}}  d \big[\widehat{ M}\big]_s\right)^{\frac{p}{2}}\right]\\
				&<+\infty.
			\end{split}
		\end{equation}
		The last term is finite due to Remark \ref{Good rmq too}.\\
		Then, taking the expectation on both sides and setting $t = 0$, $\tau = T$, and $\beta, \mu > 0$ such that $\beta, \mu > \frac{2(p-1)}{p}$, we deduce the existence of a constant $\mathfrak{c}_{\beta, \mu, p, \epsilon}$ such that
		\begin{equation}\label{for Y}
			\begin{split}
				&\mathbb{E} \int_{0}^{T}e^{\frac{p}{2}\beta A_{s} +\frac{p}{2} \mu \bar{\kappa}_{s}}\big|\widehat{Y}_s\big|^pdA_s+\mathbb{E} \int_{0}^{T}e^{\frac{p}{2}\beta A_{s} +\frac{p}{2} \mu \bar{\kappa}_{s}}\big|\widehat{Y}_s\big|^p(d\|\widehat{\kappa}\|_s+d\kappa^2_s)\\
				&+\mathbb{E}\int_{0}^{T}e^{\frac{p}{2}\beta A_{s} +\frac{p}{2} \mu \bar{\kappa}_{s}}\big|\widehat{Y}_s\big|^{p-2}\mathds{1}_{Y_s \neq 0}\big\|\widehat{Z}_s\big\|^2ds\\
				&+\mathbb{E}\int_{0}^{T}e^{\frac{p}{2}\beta A_{s} +\frac{p}{2} \mu \bar{\kappa}_{s}}\big|\widehat{Y}_s\big|^{p-2}\mathds{1}_{Y_s \neq 0}d \big[\widehat{M}\big]^c_s\\
				&+\mathbb{E}\Bigg[ \sum_{0 <s \leq T} e^{\frac{p}{2}\beta A_{s} +\frac{p}{2} \mu \bar{\kappa}_{s}} \big|\Delta \widehat{M}_s\big|^2 \left(\big|\widehat{Y}_{s-}\big|^2 \vee \big|\widehat{Y}_{s}\big|^2\right)^{\frac{p-2}{2}} \mathds{1}_{\left|\widehat{Y}_{s-}\right| \vee \left|\widehat{Y}_{s}\right| \neq 0}\Bigg] \\
				\leq & \mathfrak{c}_{\beta,\mu,p,\epsilon}\left( \mathbb{E}\left[ e^{\frac{p}{2}\beta A_{T} +\frac{p}{2} \mu \bar{\kappa}_{T}}\big|\widehat{\xi}\big|^p\right] 
				+\int_{0}^{T}e^{\beta A_{s} + \mu \bar{\kappa}_{s}} \big|\widehat{f}(s,Y^2_s,Z^2_s)\big|^pds\right. \\
				&\left. +\mathbb{E}\int_{0}^{T} e^{\beta A_{s} + \mu \bar{\kappa}_{s}} \big|g^1(s,Y^2_s)\big|^p d\|\widehat{\kappa}\|_s +\mathbb{E}\int_{0}^{T} e^{\beta A_{s} + \mu \bar{\kappa}_{s}}\big|\widehat{g}(s,Y^2_s)\big|^pd\kappa^2_s\right).
			\end{split}
		\end{equation}
		Let 
		\begin{equation*}
			\begin{split}
				\mathcal{X} &= e^{\frac{p}{2}\beta A_{T} +\frac{p}{2} \mu \bar{\kappa}_{T}} \big|\widehat{\xi}\big|^p + \int_{0}^{T} e^{\beta A_{s} + \mu \bar{\kappa}_{s}} \, \big|\widehat{f}(s, Y^2_s, Z^2_s)\big|^p \, ds \\ 
				&\qquad+ \int_{0}^{T} e^{\beta A_{s} + \mu \bar{\kappa}_{s}} \, \big|g^1(s, Y^2_s)\big| \, d\|\widehat{\kappa}\|_s+ \int_{0}^{T} e^{\beta A_{s} +\mu \bar{\kappa}_{s}} \, \big|\widehat{g}(s, Y^2_s)\big|^p \, d\kappa^2_s.
			\end{split}
		\end{equation*}
		Then, using (\ref{back2}) with $\beta, \mu > \frac{2(p-1)}{p}$, we have, a.s., for each $t \in [0, T]$
		\begin{equation}
				\begin{split}
					&e^{\frac{p}{2}\beta A_{t} +\frac{p}{2} \mu \bar{\kappa}_{t}}\big|\widehat{Y}_t\big|^p
					+\frac{c(p)}{2}\int_{t}^{T}e^{\frac{p}{2}\beta A_{s} +\frac{p}{2} \mu \bar{\kappa}_{s}}\big|\widehat{Y}_s\big|^{p-2}\mathds{1}_{Y_s \neq 0}\big\|\widehat{Z}_s\big\|^2ds\\
					&+c(p)\int_{t}^{T}e^{\frac{p}{2}\beta A_{s} +\frac{p}{2} \mu \bar{\kappa}_{s}}\big|\widehat{Y}_s\big|^{p-2}\mathds{1}_{Y_s \neq 0}d \big[\widehat{M}\big]^c_s\\
					&+c(p)\sum_{t <s \leq T} e^{\frac{p}{2}\beta A_{s} +\frac{p}{2} \mu \bar{\kappa}_{s}} \big|\Delta \widehat{M}_s\big|^2 \left(\big|\widehat{Y}_{s-}\big|^2 \vee \big|\widehat{Y}_{s-}+\Delta \widehat{M}_s\big|^2\right)^{\frac{p-2}{2}} \mathds{1}_{\left|\widehat{Y}_{s-}\right| \vee \left|\widehat{Y}_{s-}+\Delta \widehat{ M}_s\right| \neq 0}\\
					\leq &  \mathcal{X}-p\left(\Lambda_T-\Lambda_t\right)-p\left(\Xi_T-\Xi_t\right).
				\end{split}
		\end{equation}
		Using again the BDG inequality, we derive
		\begin{equation}\label{Plugg Y}
			\begin{split}
				&\mathbb{E}\left[\sup_{t \in [0,T]}e^{\frac{p}{2}\beta A_{t} +\frac{p}{2} \mu \bar{\kappa}_{t}}\big|\widehat{Y}_t\big|^p\right]
				\leq   \mathbb{E}\left[\mathcal{X}\right]+\mathfrak{c}_p\left(\mathbb{E}\left[  \Lambda\right]^{1/2}_T +\mathbb{E}\left[  \Xi\right]^{1/2}_T \right).
			\end{split}
		\end{equation}
		The term $\left[ \Lambda\right]^{1/2}_T$ can be controlled as in \cite{briand2003lp}:
		\begin{equation}\label{BDG Z}
				\begin{split}
					\mathfrak{c}_p\mathbb{E}\left[  \Lambda\right]^{1/2}_T 
					\leq & \mathfrak{c}_p  \mathbb{E}\left[\left(\sup_{t \in [0,T]}e^{\frac{p}{4}\beta A_{t} +\frac{p}{4} \mu \bar{\kappa}_{t}} \big|Y_t\big|^{\frac{p}{2}}\right)\left(\int_{0}^{T}e^{\frac{p}{2}\beta A_{s} +\frac{p}{2} \mu \bar{\kappa}_{s}} \big|\widehat{Y}_s\big|^{p-2} \mathds{1}_{Y_s\neq 0} \|\widehat{Z}_s\|^2 ds\right)^{\frac{1}{2}}\right]\\
					\leq & \frac{1}{4}\mathbb{E}\left[\sup_{t \in [0,T]}e^{\frac{p}{2}\beta A_{t} +\frac{p}{2} \mu \bar{\kappa}_{t}}\big|\widehat{Y}_t\big|^p\right]+\mathfrak{c}^2_p \mathbb{E}\int_{0}^{T}e^{\frac{p}{2}\beta A_{s} +\frac{p}{2} \mu \bar{\kappa}_{s}}\big|\widehat{Y}_s\big|^{p-2}\mathds{1}_{Y_s \neq 0}\big|\widehat{Z}_s\big|^2ds.
				\end{split}
		\end{equation}
		For the term $\left[\Xi\right]^{1/2}_T$, which is more complicated to handle, we follow \cite{kruse2016bsdes} to obtain a bound in terms of the estimation (\ref{for Y}):
		\begin{equation*}
				\begin{split}
					\mathfrak{c}_p\mathbb{E}\left[  \Xi\right]^{1/2}_T
					=&	\mathfrak{c}_p\mathbb{E}\left[  \left(\int_{0}^{T}e^{p\beta A_{s} +p \mu \bar{\kappa}_{s}} \big|\widehat{Y}_{s-}\big|^{2(p-1)} d \big[\widehat{ M}\big]_s\right)^{\frac{1}{2}}\right]\\ 
					\leq& \mathfrak{c}_p\mathbb{E}\left[  \left(\int_{0}^{T}e^{p\beta A_{s} +p \mu \bar{\kappa}_{s}} \left(\big|\widehat{Y}_{s-}\big|^{2} \vee \big|\widehat{Y}_{s}\big|^{2}\right)^{p-1} \mathds{1}_{\left|\widehat{Y}_{s-}\right| \vee \left|\widehat{Y}_{s}\right| \neq 0} d \big[\widehat{ M}\big]_s\right)^{\frac{1}{2}}\right] \\
					\leq & \mathfrak{c}_p  \mathbb{E}\left[\left(\sup_{t \in [0,T]}e^{\frac{p}{2}\beta A_{t} +\frac{p}{2} \mu \bar{\kappa}_{t}} \left(\big|\widehat{Y}_{t-}\big|^{2} \vee \big|\widehat{Y}_{t}\big|^{2}\right)^{\frac{p}{2}}\right)^{\frac{1}{2}}\right.\\
					&\qquad \left. \times \left(\int_{0}^{T}e^{\frac{p}{2}\beta A_{s} +\frac{p}{2} \mu \bar{\kappa}_{s}}\left(\big|\widehat{Y}_{s-}\big|^{2} \vee \big|\widehat{Y}_{s}\big|^{2}\right)^{\frac{p-2}{2}}  \mathds{1}_{\left|\widehat{Y}_{s-}\right| \vee \left|\widehat{Y}_{s}\right| \neq 0} d \big[\widehat{ M} \big]_s \right)^{\frac{1}{2}}\right]\\
					\leq & \frac{1}{4}\mathbb{E}\left[\sup_{t \in [0,T]}e^{\frac{p}{2}\beta A_{t} +\frac{p}{2} \mu \bar{\kappa}_{t}}\big|\widehat{Y}_t\big|^p\right]+\mathfrak{c}^2_p \mathbb{E}\int_{0}^{T}e^{\frac{p}{2}\beta A_{s} +\frac{p}{2} \mu \bar{\kappa}_{s}}\big|\widehat{Y}_s\big|^{p-2}\mathds{1}_{Y_s \neq 0}\big|\widehat{Z}_s\big|^2ds
				\end{split}
		\end{equation*}
		Using the pathwise decomposition of the bracket process $\big[\widehat{M}\big]$ (see, e.g., \cite[p. 70]{protter2005stochastic}), we have
		\begin{equation}\label{Qua}
			d\big[\widehat{M}\big]_t=d\big[\widehat{M}\big]^c_t+ \big|\Delta \widehat{M}_t\big|^2,\quad t \in [0,T],
		\end{equation}
		Therefore,
		\begin{equation*}
			\begin{split}
				&\int_{0}^{T}e^{\frac{p}{2}\beta A_{s} +\frac{p}{2} \mu \bar{\kappa}_{s}}\left(\big|\widehat{Y}_{s-}\big|^{2} \vee \big|\widehat{Y}_{s}\big|^{2}\right)^{\frac{p-2}{2}}  \mathds{1}_{\left|\widehat{Y}_{s-}\right| \vee \left|\widehat{Y}_{s}\right| \neq 0} d \big[\widehat{ M} \big]_s \\
				=&\int_{0}^{T}e^{\frac{p}{2}\beta A_{s} +\frac{p}{2} \mu \bar{\kappa}_{s}} \big|\widehat{Y}_{s}\big|^{p-2} \mathds{1}_{\left|\widehat{Y}_{s}\right|  \neq 0} d \big[\widehat{ M} \big]^c_s\\
				&+\sum_{0 < s \leq T}e^{\frac{p}{2}\beta A_{s} +\frac{p}{2} \mu \bar{\kappa}_{s}}\left(\big|\widehat{Y}_{s-}\big|^{2} \vee \big|\widehat{Y}_{s}\big|^{2}\right)^{\frac{p-2}{2}} \big|\Delta \widehat{M}_s\big|^2 \mathds{1}_{\left|\widehat{Y}_{s-}\right| \vee \left|\widehat{Y}_{s}\right| \neq 0} 
			\end{split}
		\end{equation*}
		Injecting this into the above estimation and using the basic inequality $\sqrt{ab} \leq \frac{1}{4} a + b$ for any $a, b \geq 0$, we obtain
		\begin{equation}\label{BDG M}
				\begin{split}
					\mathfrak{c}_p\mathbb{E}\left[  \Xi\right]^{1/2}_T
					\leq & \frac{1}{4}\mathbb{E}\left[\sup_{t \in [0,T]}e^{\frac{p}{2}\beta A_{t} +\frac{p}{2} \mu \bar{\kappa}_{t}}\big|\widehat{Y}_t\big|^p\right]+\mathfrak{c}_p\mathbb{E}\int_{0}^{T}e^{\frac{p}{2}\beta A_{s} +\frac{p}{2} \mu \bar{\kappa}_{s}} \big|\widehat{Y}_{s}\big|^{p-2} \mathds{1}_{\left|\widehat{Y}_{s}\right|  \neq 0} d \big[\widehat{ M} \big]^c_s\\
					&+\mathfrak{c}^2_p \mathbb{E}\sum_{0 < s \leq T}e^{\frac{p}{2}\beta A_{s} +\frac{p}{2} \mu \bar{\kappa}_{s}}\left(\big|\widehat{Y}_{s-}\big|^{2} \vee \big|\widehat{Y}_{s}\big|^{2}\right)^{\frac{p-2}{2}} \big|\Delta \widehat{M}_s\big|^2 \mathds{1}_{\left|\widehat{Y}_{s-}\right| \vee \left|\widehat{Y}_{s}\right| \neq 0} 
				\end{split}
		\end{equation}
		Plugging (\ref{BDG Z}) and (\ref{BDG M}) into (\ref{Plugg Y}), along with (\ref{for Y}), we derive
		\begin{equation}\label{Y sup}
			\begin{split}
				&\mathbb{E}\left[\sup_{t \in [0,T]}e^{\frac{p}{2}\beta A_{t} +\frac{p}{2} \mu \bar{\kappa}_{t}}\big|\widehat{Y}_t\big|^p\right]\\
				\leq & \mathfrak{c}_{\beta,\mu,p,\epsilon}\left( \mathbb{E}\left[ e^{\frac{p}{2}\beta A_{T} +\frac{p}{2} \mu \bar{\kappa}_{T}}\big|\widehat{\xi}\big|^p\right] 
				+\mathbb{E}\int_{0}^{T}e^{\beta A_{s} + \mu \bar{\kappa}_{s}}  \big|\widehat{f}(s,Y^2_s,Z^2_s)\big|^pds\right. \\
				&\left. +\mathbb{E}\int_{0}^{T} e^{\beta A_{s} +\mu \bar{\kappa}_{s}} \big|g^1(s,Y^2_s)\big|^pd\|\widehat{\kappa}\|_s  +\mathbb{E}\int_{0}^{T} e^{\beta A_{s} + \mu \bar{\kappa}_{s}}\big|\widehat{g}(s,Y^2_s)\big|^pd\kappa^2_s\right).
			\end{split}
		\end{equation}
		
		Showing the estimation for the remaining term 
		$$
		\mathbb{E}\left[\left(\int_{0}^{T} e^{\beta A_{s} +\mu \bar{\kappa}_{s}} \, \big\|\widehat{Z}_s\big\|^2 \, ds\right)^{\frac{p}{2}}\right] 
		+ \mathbb{E}\left[\left(\int_{0}^{T} e^{\beta A_{s} + \mu \bar{\kappa}_{s}} \, d\big[\widehat{M}\big]_s\right)^{\frac{p}{2}}\right]
		$$
		is the last step in proving the assertion. To this end, since $p \in (1, 2)$, we apply the equality $
		\mathds{1}_{\widehat{Y}_s=0} \left\{ \big\|\widehat{Z}_s\big\|^2 \, ds + d\big[\widehat{M}\big]_s \right\} = 0
		$ on $[0, T]$ (we refer to \cite[Lemma 8]{kruse2016bsdes} for a detailed proof). Using this result, along with Young's inequality, we obtain
		\begin{equation*}
				\begin{split}
					&\mathbb{E}\left[\left(\int_{0}^{T}e^{\beta A_{s} +\mu \bar{\kappa}_{s}}\big\|\widehat{Z}_s\big\|^2 ds\right)^{\frac{p}{2}}\right]\\
					=&\mathbb{E}\left[\left(\int_{0}^{T}e^{\beta A_{s} + \mu \bar{\kappa}_{s}} \big\|\widehat{Z}_s\big\|^2 \mathds{1}_{\{\widehat{Y}_s\neq0\}} ds\right)^{\frac{p}{2}}\right]\\
					=&\mathbb{E}\left[\left(\int_{0}^{T}e^{\frac{2-p}{2}\beta A_{s} +\frac{2-p}{2} \mu \bar{\kappa}_{s}} \big|\widehat{Y}_s\big|^{2-p} \big(e^{\frac{p}{2}\beta A_{s} +\frac{p}{2} \mu \bar{\kappa}_{s}}\big\|\widehat{Z}_s\big\|^2 \big|\widehat{Y}_s\big|^{p-2} \mathds{1}_{\{\widehat{Y}_s\neq0\}}\big) ds\right)^{\frac{p}{2}}\right]\\
					\leq& \mathbb{E}\left[\left(\sup_{t \in [0,T]}e^{\frac{p(2-P)}{2}\beta A_{t} +\frac{p(2-p)}{2} \mu \bar{\kappa}_{t}} \big|\widehat{Y}_t\big|^{\frac{p(2-p)}{2}}\right)\right. \\
					&\left. \quad ~\times\left(\int_{0}^{T}e^{\frac{p}{2}\beta A_{s} +\frac{p}{2} \mu \bar{\kappa}_{s}} \big\|\widehat{Z}_s\big\|^2 \big|\widehat{Y}_s\big|^{p-2} \mathds{1}_{\{\widehat{Y}_s\neq0\}} ds\right)^{\frac{p}{2}}\right]\\
					\leq &\frac{2-p}{2}\mathbb{E}\left[\sup_{t \in [0,T]}e^{\frac{p}{2}\beta A_{t} +\frac{p}{2} \mu \bar{\kappa}_{t}}\big|\widehat{Y}_{t}\big|^p\right]+\frac{p}{2} \mathbb{E}\left[\int_{0}^{T}e^{\frac{p}{2}\beta A_{s} +\frac{p}{2} \mu \bar{\kappa}_{s}}\big|\widehat{Y}_{s}\big|^{p-2} \big\|\widehat{Z}_s\big\|^2 \mathds{1}_{\{\widehat{Y}_{s}\neq 0\} }ds\right].
				\end{split}
		\end{equation*}
		Similarly, we can show that
		\begin{equation*}
				\begin{split}
					&\mathbb{E}\left[\left(\int_{0}^{T}e^{\beta A_{s} + \mu \bar{\kappa}_{s}}d\big[\widehat{M}\big]^c_s \right)^{\frac{p}{2}}\right]\\
					\leq &\frac{2-p}{2}\mathbb{E}\left[\sup_{t \in [0,T]}e^{\frac{p}{2}\beta A_{t} +\frac{p}{2} \mu \bar{\kappa}_{t}}\big|\widehat{Y}_{t}\big|^p\right]+\frac{p}{2} \mathbb{E}\left[\int_{0}^{T}e^{\frac{p}{2}\beta A_{s} +\frac{p}{2} \mu \bar{\kappa}_{s}}\big|\widehat{Y}_{s}\big|^{p-2}  \mathds{1}_{\big|\widehat{Y}_{s}\big|\neq 0 }d\big[\widehat{M}\big]^c_s\right].
				\end{split}
		\end{equation*}
		By virtue of (\ref{Qua}), it remains to show a similar estimation for the quadratic jump part of the state process $\widehat{Y}$ given by $\big| \Delta \widehat{M} \big|^2$. To this end, we will use an approximating procedure via a smooth function that is widely considered in the literature (see, e.g., \cite[Lemma 2.2]{briand2003lp}, \cite[Lemma 7]{kruse2016bsdes}, or \cite[Lemma 2.2]{ElJamali}).\\ Let $\varepsilon > 0$, and consider the function $\nu_\varepsilon : \mathbb{R} \to \mathbb{R}_+$ defined by $\nu_\varepsilon(y) = \sqrt{\left|y\right|^2 + \varepsilon^2}$. Then, for any $q > 0$, we have
		\begin{equation*}
			\begin{split}
				e^{\frac{2-p}{2}\beta A_{s} +\frac{2-p}{2} \mu \bar{\kappa}_{s}} \left(\nu_\varepsilon(y)\right)^q&=\left(\left(e^{\frac{2-p}{q}\beta A_{s} +\frac{2-p}{q} \mu \bar{\kappa}_{s}}\left(\big|y\big|^2+\varepsilon^2\right)\right)^{\frac{1}{2}}\right)^q\\ &\leq\left(\nu_{\widehat{\varepsilon}^{p,q}_s}\left(e^{\frac{2-p}{2q}\beta A_{s} +\frac{2-p}{2q} \mu \bar{\kappa}_{s}}y\right)\right)^q,
			\end{split}
		\end{equation*}
		where $\widehat{\varepsilon}_{p,q} = \varepsilon \mathcal{C}^{\ast}_{p,q}$ with $\mathcal{C}^{\ast}_{p,q} := \esssup_{\omega \in \Omega} e^{\frac{p}{2}\beta A_{T}(\omega) +\frac{p}{2} \mu \bar{\kappa}_{T}(\omega)}$. Without loss of generality, instead of replacing $A$ with $A \wedge k$ and then passing to the limit using the monotone convergence theorem, we may assume that $A$ is bounded. Therefore, we have $\mathcal{C}^{\ast}_{p,q} < +\infty$, and consequently, $\lim\limits_{\varepsilon \downarrow 0} \widehat{\varepsilon}_{p,q} = 0$. Furthermore, 
		$$
		\lim\limits_{\varepsilon \downarrow 0} \nu_{\widehat{\varepsilon}_{p,q}} \left( e^{\frac{2-p}{2q}\beta A_{s} +\frac{2-p}{2q} \mu \bar{\kappa}_{s}} \, y \right) = e^{\frac{2-p}{2q}\beta A_{s} +\frac{2-p}{2q} \mu \bar{\kappa}_{s}} \, |y| \, \mathds{1}_{y \neq 0} \text{ a.s., } q > 0.$$
		To simplify notation, we denote $\mathcal{X}_\ast = \sup_{t \in [0,T]} \big| \mathcal{X}_t \big|$ for any RCLL process $\mathcal{X} = (\mathcal{X}_t)_{t \leq T}$. Then, using the H\"older's and Young's inequalities 
		$$
		\mathbb{E} \left[ A^{p(2-p)/2} \, B^{p/2} \right] \leq \left( \mathbb{E} \left[ A^p \right] \right)^{(2-p)/2} \left( \mathbb{E} \left[ B \right] \right)^{p/2} \leq \frac{2-p}{2} \, \mathbb{E} \left[ A^p \right] + \frac{p}{2} \, \mathbb{E} \left[ B \right]
		$$ 
		for some random variables $A, B \geq 0$, we have
		\begin{equation}\label{LbMC}
				\begin{split}
					&\mathbb{E}\left[\left(\sum_{0 < s \leq T} e^{\beta A_{s} + \mu \bar{\kappa}_{s}} \big|\Delta \widehat{M}_s\big|^2 \right)^{\frac{p}{2}}\right]\\
					= &\mathbb{E}\left[\left(\sum_{0 < s \leq T} e^{\beta A_{s} +\mu \bar{\kappa}_{s}} \left( \nu_\varepsilon\left(\big|Y_{s-}\big| \vee \big|Y_{s}\big| \right)\right)^{2-p}\left( \nu_\varepsilon\left(\big|Y_{s-}\big| \vee \big|Y_{s}\big| \right)\right) ^{p-2} \big|\Delta \widehat{M}_s\big|^2 \right)^{\frac{p}{2}}\right]\\
					&\leq \mathbb{E}\left[\left(\nu_{\widehat{\varepsilon}_{p,2-p}}\left(\left(e^{\frac{\beta}{2} A +\frac{\mu}{2} \bar{\kappa}} \widehat{Y}\right)_\ast\right)\right)^{\frac{p(2-p)}{2}}\right. \\
					&\left. \qquad\qquad \times\left(\sum_{0 <s \leq T}e^{\frac{p}{2}\beta A_{s} +\frac{p}{2} \mu \bar{\kappa}_{s}} \left( \nu_\varepsilon\left(\big|Y_{s-}\big| \vee \big|Y_{s}\big| \right)\right)^{2-p} \big|\Delta \widehat{M}_s\big|^2 \right)^{\frac{p}{2}}\right]\\
					&\leq \left(\mathbb{E}\left[\left(\nu_{\widehat{\varepsilon}^{p,2-p}_s}\left(\left(e^{\frac{\beta}{2} A +\frac{\mu}{2} \bar{\kappa}} \widehat{Y}\right)_\ast\right)\right)^{p}\right]\right)^{\frac{2-p}{2}}\\ 
					&\qquad\qquad \times\left(\mathbb{E}\left[\sum_{0 <s \leq T}e^{\frac{p}{2}\beta A_{s} +\frac{p}{2} \mu \bar{\kappa}_{s}} \left( \nu_\varepsilon\left(\big|Y_{s-}\big| \vee \big|Y_{s}\big| \right)\right)^{p-2} \big|\Delta \widehat{M}_s\big|^2\right] \right)^{\frac{p}{2}}\\
					& \leq \frac{2-p}{2} \mathbb{E}\left[\left(\nu_{\widehat{\varepsilon}^{p,2-p}_s}\left(\left(e^{\frac{\beta}{2} A +\frac{\mu}{2} \bar{\kappa}} \widehat{Y}\right)_\ast\right)\right)^{p}\right]\\
					&\qquad\qquad +\frac{p}{2} \mathbb{E}\left[\sum_{0 <s \leq T}e^{\frac{p}{2}\beta A_{s} +\frac{p}{2} \mu \bar{\kappa}_{s}} \left( \nu_\varepsilon\left(\big|Y_{s-}\big| \vee \big|Y_{s}\big| \right)\right)^{p-2} \big|\Delta \widehat{M}_s\big|^2\right]
				\end{split}
		\end{equation}
		We know that 
		$$\lim\limits_{\varepsilon \downarrow 0}\left(\nu_{\widehat{\varepsilon}^{p,2-p}_s}\left(\left(e^{\frac{\beta}{2} A +\frac{\mu}{2}  \bar{\kappa}} \widehat{Y}\right)_\ast\right)\right)^{p}=\sup_{t \in [0,T]}e^{\frac{p}{2}\beta A_{t} +\frac{p}{2} \mu \bar{\kappa}_{t}} \big|\widehat{Y}_t\big|^p \mbox{ a.s.}$$  
		and that 
		$$\lim\limits_{\varepsilon \downarrow 0}\left(\nu_\varepsilon\left(\big|\widehat{Y}_{s-}\big| \vee \big|\widehat{Y}_{s}\big| \right)\right)^{p-2} \nearrow \left(\big|Y_{s-}\big| \vee \big|Y_{s}\big|\right)^{p-2}\mathds{1}_{\big|\widehat{Y}_{s-}\big| \vee \big|\widehat{Y}_{s}\big| \neq 0} \mbox{ a.s.,}$$ 
		since $p < 2$.\\ 
		Letting $\varepsilon \to 0$ and then applying the Lebesgue dominated convergence theorem for the left term in the last inequality of the estimation (\ref{LbMC}) and the monotone convergence theorem for the right term in the last inequality of the estimation (\ref{LbMC}) (see also the proof of Lemma 9 in \cite{kruse2016bsdes}), we obtain
		\begin{equation*}
			\begin{split}
				&\mathbb{E}\left[\left(\sum_{0 < s \leq T} e^{\beta A_{s} +\mu \bar{\kappa}_{s}} \big|\Delta \widehat{M}_s\big|^2 \right)^{\frac{p}{2}}\right]\\
				\leq &\frac{2-p}{2}\mathbb{E}\left[\sup_{t \in [0,T]}e^{\frac{p}{2}\beta A_{t} +\frac{p}{2} \mu \bar{\kappa}_{t}}\big|\widehat{Y}_{t}\big|^p\right]\\
				&\qquad+\frac{p}{2} \mathbb{E}\left[\sum_{0 <s \leq T}e^{\frac{p}{2}\beta A_{s} +\frac{p}{2} \mu \bar{\kappa}_{s}}\left(\big|Y_{s-}\big|^2 \vee \big|Y_{s}\big|^2\right)^{\frac{p-2}{2}}\mathds{1}_{\big|\widehat{Y}_{s-}\big| \vee \big|\widehat{Y}_{s}\big| \neq 0} \big|\Delta \widehat{M}_s\big|^2\right].
			\end{split}
		\end{equation*}
		Finally, the proof is completed by applying the inequalities (\ref{for Y}) and (\ref{Y sup}).
	\end{proof}

	From Proposition \ref{propo1}, we obtain the uniqueness of the solution.
	\begin{corollary}\label{Uniquen}
		Let $(\xi, f, g, \kappa)$ be any set of data satisfying assumption \textbf{(H-M)$_\mathbf{p}$}. Then, there exists at most one triplet of processes $(Y_t, Z_t, M_t)_{t \leq T}$ corresponding to the $\mathbb{L}^p$-solution of the GBSDE (\ref{basic GBSDE}) associated with $(\xi, f, g, \kappa)$.
	\end{corollary}

	\section{Existence}
	\label{sec4}
	To establish the existence and uniqueness of $\mathbb{L}^p$-solutions for $p \in (1,2)$ of the GBSDE (\ref{basic GBSDE}), we first need to prove the existence and uniqueness of $\mathbb{L}^2$-solutions. While this result may have appeared in previous works (see, e.g., \cite{ElmansouriElOtmani,elmansourielotmani2024}), we could not find an equivalent result for the GBSDE (\ref{basic GBSDE}) under our general condition \textbf{(H-M)$_\mathbf{2}$}, which assumes stochastic monotonicity, Lipschitz continuity, and linear growth in a general filtration. To address this gap, we begin by analyzing the $\mathbb{L}^2$ case.
	
	\subsection{Existence of $\mathbb{L}^2$-solutions}
	It is evident that in the $\mathbb{L}^2$ case, the situation is more regular than in the $\mathbb{L}^p$ case for $p \in (1,2)$. The classical multidimensional It\^o's formula (see, e.g., Theorem 33 in \cite[pp. 81-82]{protter2005stochastic}) can be directly applied to the function $\mathbb{R}^d \ni y \mapsto |y|^2$, and the jump part of the process $Y$ (i.e., $\Delta M$) can be more easily controlled. With this in mind, we only provide the key elements of the proof.
	
	\subsubsection{A priori estimates and uniqueness}
	Let $(Y^1_t, Z^1_t, M^1_t)_{t \leq T}$ and $(Y^2_t, Z^2_t, M^2_t)_{t \leq T}$ be two $\mathbb{L}^2$-solutions of the GBSDE (\ref{basic GBSDE}) associated with the data $(\xi^1, f^1, g^1, \kappa^1)$ and $(\xi^2, f^2, g^2, \kappa^2)$, respectively, satisfying condition \textbf{(H-M)$_\mathbf{2}$}. Define $\mathcal{R} = \mathcal{R}^1 - \mathcal{R}^2$ with $\mathcal{R} \in \left\{Y, Z, M, \xi, f, \kappa \right\}$. Then, we have
	\begin{proposition}\label{propo2}
		For any $\beta, \mu > 0$, there exists a constant $\mathfrak{c}_{\beta, \mu}$ such that
		\begin{equation*}
			\begin{split}
				&\mathbb{E}\left[\sup_{t \in [0,T]}e^{\beta A_{t} + \mu \bar{\kappa}_{t}}\big|\widehat{Y}_t\big|^2\right]+\mathbb{E} \int_{0}^{T}e^{\beta A_{s} + \mu \bar{\kappa}_{s}}\big|\widehat{Y}_s\big|^2 (dA_s+d\|\widehat{\kappa}\|_s+d\kappa^2_s)\\
				&+\mathbb{E} \int_{0}^{T}e^{\beta A_{s} + \mu \bar{\kappa}_{s}}\big\|\widehat{Z}_s\big\|^2ds  
				+\mathbb{E} \int_{0}^{T}e^{\beta A_{s} + \mu \bar{\kappa}_{s}}d\big[\widehat{M}\big]_s \\
				\leq & \mathfrak{c}_{\beta,\mu}\left( \mathbb{E}\left[ e^{\beta A_{T} + \mu \bar{\kappa}_{T}}\big|\widehat{\xi}\big|^2\right] 
				+\int_{0}^{T}e^{\beta A_{s} + \mu \bar{\kappa}_{s}} \left| \frac{\widehat{f}(s,Y^2_s,Z^2_s)}{a_s}\right|^2 ds\right. \\
				&\left.  \qquad+\mathbb{E}\int_{0}^{T} e^{\beta A_{s} + \mu \bar{\kappa}_{s}} \big|g^1(s,Y^2_s)\big|^2 d\|\widehat{\kappa}\|_s+\mathbb{E}\int_{0}^{T} e^{\beta A_{s} + \mu \bar{\kappa}_{s}}\big|\widehat{g}(s,Y^2_s)\big|^2 d\kappa^2_s\right).
			\end{split}
		\end{equation*}
	\end{proposition}
	\begin{proof}
		By applying It\^o's formula to the special semimartingale $e^{\beta A+ \mu \bar{\kappa}} \big|\widehat{Y}\big|^2$, where $\widehat{Y}$ satisfies the dynamic (\ref{Khdj}), we have
		\begin{equation}\label{Ito L2}
			\begin{split}
				&e^{\beta A_{t} + \mu \bar{\kappa}_{t}}\big|\widehat{Y}_{t}\big|^2+\beta \int_{t }^{T}e^{\beta A_{s} + \mu \bar{\kappa}_{s}}\big|\widehat{Y}_s\big|^2dA_s+\mu  \int_{t  }^{T}e^{\beta A_{s} + \mu \bar{\kappa}_{s}}\big|\widehat{Y}_s\big|^2 (d\|\widehat{\kappa}\|_s+d\kappa^2_s)\\
				&+\int_{t}^{T}e^{\beta A_{s} + \mu \bar{\kappa}_{s}}\big\|\widehat{Z}_s\big\|^2ds+\int_{t }^{T}e^{\beta A_{s} + \mu \bar{\kappa}_{s}} d \big[\widehat{M}\big]_s\\
				= &e^{\beta A_{T} + \mu \bar{\kappa}_{T}}\big|\widehat{\xi}\big|^2
				+2\int_{t}^{T}e^{\beta A_{s} + \mu \bar{\kappa}_{s}} \widehat{Y}_s \left(f^1(s,Y^1_s,Z^1_s)-f^2(s,Y^2_s,Z^2_s)\right) ds\\
				&+2\int_{t}^{T} e^{\beta A_{s} + \mu \bar{\kappa}_{s}} \widehat{Y}_s \left(g^1(s,Y^1_s)d\kappa^1_s-g^2(s,Y^2_s)d\kappa^2_s\right)\\
				&-2\int_{t}^{T}e^{\beta A_{s} + \mu \bar{\kappa}_{s}} \widehat{Y}_s  \widehat{Z}_s dW_s-2\int_{t}^{T}e^{\beta A_{s} + \mu \bar{\kappa}_{s}} \widehat{Y}_{s-}  d\widehat{M}_s.
			\end{split}
		\end{equation}
		From the assumptions (H2)–(H4) on $f$ and $g$, Remark \ref{rmq essential}, and the basic inequality $2ab \leq \varepsilon a^2 + \frac{1}{\varepsilon} b^2$ for any $\varepsilon > 0$, we have
		\begin{equation*}
			\begin{split}
				&2\widehat{Y}_s \left(f^1(s,Y^1_s,Z^1_s)-f^2(s,Y^2_s,Z^2_s)\right) ds\\
				&\leq \frac{\beta}{2}\big| \widehat{Y}_s \big|^2 dA_s+ \frac{1}{2} \big\| \widehat{Z}_s \big\|^2ds+\frac{2}{\beta}\left| \frac{\widehat{f}(s,Y^2_s,Z^2_s)}{a_s}\right|^2 ds
			\end{split}
		\end{equation*}
		and
		\begin{equation*}
			\begin{split}
				&2\widehat{Y}_s\left(g^1(s,Y^1_s)d\kappa^1_s-g^2(s,Y^2_s)d\kappa^2_s\right)\\
				\leq& 2\big|\widehat{Y}_s\big| \big|g^1(s,Y^2_s)\big|d\|\widehat{\kappa}\|_s+2\big|\widehat{Y}_s\big| \big|g^1(s,Y^2_s)-g^2(s,Y^2_s)\big|d\kappa^2_s\\
				\leq & \frac{\mu}{2}\big|\widehat{Y}_s\big|^2\left(d\|\widehat{\kappa}\|_s+d\kappa^2_s\right)+\frac{2}{\mu}\left(\big|g^1(s,Y^2_s)\big|^2 d\|\widehat{\kappa}\|_s+\big|g^1(s,Y^2_s)-g^2(s,Y^2_s)\big|^2 d\kappa^2_s\right)
			\end{split}
		\end{equation*}
		Additionally, by following a similar approach to the estimations (\ref{L2})–(\ref{L3}), we can prove that the local martingale
		$$
		\int_{0}^{\cdot}e^{\beta A_{s} + \mu \bar{\kappa}_{s}} \widehat{Y}_s  \widehat{Z}_s dW_s+\int_{0}^{\cdot}e^{\beta A_{s} + \mu \bar{\kappa}_{s}} \widehat{Y}_{s-}  d\widehat{M}_s
		$$
		is a uniformly integrable martingale with zero expectation due to the BDG inequality for $p=2$ (recall that the sum of two local martingales remains a local martingale; see, e.g., \cite[Theorem 48, p. 37]{protter2005stochastic}). Moreover, note that in the formula (\ref{Ito L2}), there is no jump term $\sum \Delta \widehat{M}$ involved in the quadratic variation, unlike in (\ref{Coming back to this Itos formula}). Consequently, the expression (\ref{Ito L2}) is simpler to handle, as there is no need to control each martingale part separately.\\
		Next, using the above estimations on $f$ and $g$, along with the martingale property and following similar arguments to those used in the proof of Proposition \ref{propo1} (see also \cite[Proposition 1]{elmansourielotmani2024}), we derive the desired result.
	\end{proof}

	Following Corollaries \ref{Uniquen} and \ref{coro Lp}, and using Proposition \ref{propo2}, we derive:
	\begin{corollary}\label{coro L2}
		Let $(\xi, f, g, \kappa)$ be any set of data satisfying assumption \textbf{(H-M)$_\mathbf{2}$}. Then, there exists at most one triplet of processes $(Y_t, Z_t, M_t)_{t \leq T}$ corresponding to the $\mathbb{L}^2$-solution of the GBSDE (\ref{basic GBSDE}) associated with $(\xi, f, g, \kappa)$. Moreover, for any $\beta, \mu > 0$, there exists a constant $\mathfrak{c}_{\beta, \mu}$ such that
		\begin{equation*}
				\begin{split}
					&\mathbb{E}\left[\sup_{t \in [0,T]}e^{\beta A_{t} + \mu {\kappa}_{t}}\big|{Y}_t\big|^2\right]+\mathbb{E} \int_{0}^{T}e^{\beta A_{s} + \mu {\kappa}_{s}}\big|{Y}_s\big|^2 dA_s+\mathbb{E} \int_{0}^{T}e^{\beta A_{s} + \mu {\kappa}_{s}}\big|{Y}_s\big|^2 d\kappa_s\\
					&+\mathbb{E}\int_{0}^{T}e^{\beta A_{s} + \mu {\kappa}_{s}}\big\|{Z}_s\big\|^2ds+\mathbb{E}\int_{0}^{T}e^{\beta A_{s} + \mu {\kappa}_{s}} d\big[{M}\big]_s \\
					\leq & \mathfrak{c}_{\beta,\mu}\left( \mathbb{E}\left[ e^{\beta A_{T} + \mu {\kappa}_{T}}\big|{\xi}\big|^2\right] 
					+\int_{0}^{T}e^{\beta A_{s} + \mu {\kappa}_{s}} \left|\frac{\varphi_s}{a_s}\right|^2 ds+\mathbb{E}\int_{0}^{T} e^{\beta A_{s} + \mu {\kappa}_{s}}\big|\psi_s\big|^2 d\kappa_s\right).
				\end{split}
		\end{equation*}
	\end{corollary}
	\subsubsection{Existence}
	In this part, we establish the existence of an $\mathbb{L}^2$-solution for the GBSDE (\ref{basic GBSDE}). To achieve this, we divide the proof into two steps:
	\begin{enumerate}
		\item In the first step, we consider a generator that is independent of the $z$-variable and depends only on $y$. Using the Yosida approximation method, we derive the existence and uniqueness of an $\mathbb{L}^2$-solution.
		
		\item In the second step, we address the case of a global generator $f$ that depends on both $(y, z)$. By applying a fixed-point theorem in an appropriate Banach space, we establish the existence of an $\mathbb{L}^2$-solution.
	\end{enumerate}

	\paragraph{The case when the coefficient $f$ is independent of $z$}
	Set $\mathfrak{f}(t,y) := f(t,y,z)$ a.s. for any $(t, y, z) \in [0,T] \times \mathbb{R}^d \times \mathbb{R}^{d \times k}$. We consider the special case of the GBSDE (\ref{basic GBSDE}) associated with $(\xi, \mathfrak{f}, g, \kappa)$, which is expressed as follows:
	\begin{equation}\label{special GBSDE}
		Y_t=\xi+\int_{t}^{T}\mathfrak{f}(s,Y_s)ds+\int_{t}^{T}g(s,Y_s)d\kappa_s-\int_{t}^{T}Z_s dW_s-\int_{t}^{T}dM_s,~t \in [0,T].
	\end{equation}
	Then, we have the following lemma:
	\begin{lemma}\label{Special GRBSDE lem}
		If we assume that $\mathfrak{f}$ and $g$ are uniformly Lipschitz with respect to $y$, then the GBSDE (\ref{special GBSDE}) admits a unique solution.
	\end{lemma}
	\begin{proof}
		The proof of the existence of a unique solution $\big(Y_t, Z_t, M_t\big)_{t \leq T}$ in $\mathcal{E}^2_{0, \mu}$ for the GBSDE (\ref{special GBSDE}) is provided in \cite[Theorem 4.1]{ElmansouriElOtmani} (see also Proposition 4.2 in the same reference). Additionally, by considering the sequence of stopping times $\left\{\tau_k\right\}_{k \geq 1}$ defined by $\tau_k := \inf\left\{t \geq 0 : A_t \geq k\right\} \wedge T$, which is clearly an increasing sequence converging to $T$ a.s., we can apply the result of Corollary \ref{coro L2} on each interval $[0, \tau_k]$. By passing to the limit as $k \to +\infty$ and using Fatou's Lemma along with assumptions (H5)–(H6), we conclude that $(Y_t, Z_t, M_t)_{t \leq T}$ is the unique $\mathbb{L}^2$-solution of the GBSDE (\ref{special GBSDE}).
	\end{proof}
	
	Now, we aim to establish the existence and uniqueness of a solution for the GBSDE (\ref{special GBSDE}) under more general assumptions on the data $(\xi, \mathfrak{f}, g, \kappa)$. Specifically, we have the following proposition:
	\begin{proposition}\label{propo3}
		Suppose that $(\xi, \mathfrak{f}, g, \kappa)$ satisfies assumption \textbf{(H-M)$_\mathbf{2}$}. Then, the GBSDE (\ref{special GBSDE}) admits a unique $\mathbb{L}^2$-solution.
	\end{proposition}

	\begin{proof}
		From Remark \ref{rmq essential}, the assumptions (H2)–(H3) can be reformulated as follows: For any $t \in [0,T]$, $y, y^\prime \in \mathbb{R}^d$, and $z \in \mathbb{R}^{d \times k}$, we have:
		\begin{equation*}
			\left\lbrace 
			\begin{split}
				\left(y-y^\prime\right)\left(f(t,y,z)-f(t,y^\prime,z)\right) &\leq 0,\\
				\left(y-y^\prime\right)\left(g(t,y)-g(t,y^\prime)\right)&\leq 0.
			\end{split}
			\right. 
		\end{equation*}
		Then, using assumption (H7), we deduce that the functions $y \mapsto -\mathfrak{f}(t,y): \mathbb{R}^d \to \mathbb{R}^d$ and $y \mapsto -g(t,y): \mathbb{R}^d \to \mathbb{R}^d$ are monotone and continuous. It follows (see \cite[p. 524]{pardoux2014stochastic}) that, for every $(\omega, t, y) \in \Omega \times [0,T] \times \mathbb{R}$ and $\delta > 0$, there exist unique functions $\mathfrak{f}_\delta = \mathfrak{f}_\delta(\omega, t, y) \in \mathbb{R}^d$ and $g_\delta = g_\delta(\omega, t, y) \in \mathbb{R}^d$ such that
		\begin{equation}\label{fg}
			\mathfrak{f}(\omega,t,y+\delta \mathfrak{f}_\delta)=\mathfrak{f}_\delta,\qquad \text{and}\qquad
			g(\omega,t,y+\delta g_\delta)=g_\delta.
		\end{equation}
		From \cite[Proposition 6.7, Annex B]{pardoux2014stochastic}, the functions $\mathfrak{f}_\delta(\cdot, \cdot, y) : \Omega \times [0,T] \to \mathbb{R}^d$ and $g_\delta(\cdot, \cdot, y) : \Omega \times [0,T] \to \mathbb{R}^d$ are progressively measurable stochastic processes. Moreover, for all $\delta > 0$, $t \in [0,T]$, and $y, y' \in \mathbb{R}^d$, they satisfy the following properties a.s.:
		\begin{equation}\label{properties of Yosida approx}
			\left\{
			\begin{array}{ll}
				\displaystyle(y-y')\big(\mathfrak{f}_\delta(t,y)-\mathfrak{f}_\delta(t,y')\big)\leq 0,  \\
				\displaystyle|\mathfrak{f}_\delta(t,y)-\mathfrak{f}_\delta(t,y')|\leq \frac{2}{\delta}|y-y'|,\\
				\left|\mathfrak{f}_{\delta}\left(t,y\right) \right| \leq \left|\mathfrak{f}\left(t,y\right) \right|\\
			\end{array}
			\right. \mbox{ and }
			\left\{
			\begin{array}{ll}
				\displaystyle (y-y')\big(g_\delta(t,y)-g_\delta(t,y')\big)\leq 0, \\
				\displaystyle |g_\delta(t,y)-g_\delta(t,y')|\leq \frac{2}{\delta}|y-y'|,\\
				\left|g_{\delta}\left(t,y\right) \right| \leq \left|g\left(t,y\right) \right|
			\end{array}
			\right.
		\end{equation}
		Moreover, the coefficient $\mathfrak{f}_{\delta}$ and $g_{\delta}$ satisfies
		\begin{equation}
			\left\{
			\begin{split}
				&\left(y-y^{\prime}\right)\left(\mathfrak{f}_{\delta}\left(t,y\right)-\mathfrak{f}_{\delta'}\left(t,y'\right)\right)
				\leq\left( \delta+\delta'\right)\mathfrak{f}_{\delta}\left(t,y\right) \mathfrak{f}_{\delta'}\left(t,y'\right)\\
				&\left(y-y^{\prime}\right)\left(g_{\delta}\left(t,y\right)-g_{\delta'}\left(t,y'\right)\right)
				\leq\left(\delta+\delta'\right) g_{\delta}\left(t,y\right) g_{\delta'}\left(t,y'\right)
			\end{split}
			\right.
			\label{approximating equation}
		\end{equation}
		
		Let $\delta \in (0,1]$. From (\ref{properties of Yosida approx}), the functions $y \mapsto \mathfrak{f}_{\delta}(t,y)$ and $y \mapsto g_{\delta}(t,y)$ are uniformly Lipschitz with the same constant $\frac{2}{\delta}$. By Lemma \ref{Special GRBSDE lem}, we deduce that the approximating equation
		\begin{equation}\label{Yosida approximation equation}
			Y^{\delta}_{t}=\xi +\int_{t}^{T} \mathfrak{f}_{\delta}(s,Y^{\delta}_{s})ds+\int_{t}^{T} g_{\delta}(s,Y^{\delta}_{s})d\kappa_s
			-\int_{t}^{T}Z^{\delta}_{s}dW_{s}-\int_{t}^{T}dM^{\delta}_s,~t \in [0,T].
		\end{equation}
		has a unique $\mathbb{L}^2$-solution $(Y^{\delta},Z^{\delta},M^{\delta})$.
		
		Now, let $\left\{Y^{\delta}, Z^{\delta}, M^{\delta}\right\}_{\delta \in (0,1]}$ be a family of solutions to the GBSDE (\ref{Yosida approximation equation}) associated with $\left(\xi, \mathfrak{f}_{\delta}, g_{\delta}, \kappa\right)$. Based on the properties satisfied by the Yosida approximating drivers $\left\{\mathfrak{f}_{\delta}, g_{\delta}\right\}_{\delta \in (0,1]}$ given in (\ref{properties of Yosida approx}), we can conclude from Corollary \ref{coro L2} that $\left\{Y^{\delta}, Z^{\delta}, M^{\delta}\right\}_{\delta \in (0,1]}$ satisfies
		\begin{equation}
				\begin{split}
					&\mathbb{E}\left[ \sup_{t \in  [0,T]} e^{\beta A_{t} + \mu {\kappa}_{t}}|Y^{\delta}_t|^2\right] +\mathbb{E}\int_{0}^{T} e^{\beta A_{s} + \mu {\kappa}_{s}}|Y_s|^2 dA_s + \mathbb{E}	\int_{0}^{T} e^{\beta A_{s} + \mu {\kappa}_{s}}|Y^{\delta}_s|^2 d\kappa_s   \\
					&\qquad+ \mathbb{E}\int_{0}^{T}e^{\beta A_{s} + \mu {\kappa}_{s}}\left\|Z^{\delta}_s \right\|^2  ds+\mathbb{E}\int_{0}^{T}e^{\beta A_{s} + \mu {\kappa}_{s}} d\left[ M^{\delta}\right]_s 
					\\
					&\leq \mathfrak{c}_{\beta,\mu} \left(\mathbb{E} \left[ e^{\beta A_{T} + \mu {\kappa}_{T}} |\xi|^2\right] +\mathbb{E}\int_{0}^{T}e^{\beta A_{s} + \mu {\kappa}_{s}} \left|\dfrac{\varphi_s}{a_s} \right|^2ds+\mathbb{E}\int_{0}^{T}e^{\beta A_{s} + \mu {\kappa}_{s}}\left|\psi_s\right|^2d\kappa_s \right) .
				\end{split}
			\label{Uniform estimation}
		\end{equation}
		On the other hand, using H\"older's inequality, the uniform estimation (\ref{Uniform estimation}), the last property in (\ref{properties of Yosida approx}), and assumption \textsc{(H5)}, we can write
		\begin{equation}
				\begin{split}
					&\mathbb{E}\int_{0}^{T}e^{\beta A_{s} + \mu {\kappa}_{s}}\left| \mathfrak{f}_{\delta}(s,Y^{\delta}_s) \right|^2 ds
					\leq2 \mathbb{E}\int_{0}^{T}e^{\beta A_{s} + \mu {\kappa}_{s}}\left\{\left|\varphi_s \right|^2ds+\left|Y^{\delta}_s\right|^2 dA_s\right\}\leq \mathfrak{C}_{\beta,\mu}.
				\end{split}
			\label{estimation for hdelta}
		\end{equation}
		Similarly, we can derive
		\begin{equation}
				\begin{split}
					&\mathbb{E}\int_{0}^{T}e^{\beta A_{s} + \mu {\kappa}_{s}}\left| g_{\delta}(s,Y^{\delta}_s) \right|^2 d\kappa_s \leq  2\mathbb{E}\int_{0}^{T}e^{\beta A_{s} + \mu {\kappa}_{s}}\left\{\psi^2_s +\Gamma^2\left| Y^{\delta}_s \right|^2   \right\}d\kappa_s\leq \mathfrak{C}_{\beta,\mu}.
				\end{split}
			\label{estimation for gdelta}
		\end{equation}
		Here, the constant $\mathfrak{C}_{\beta,\mu}$ is independent of the parameter $\delta$. 
		
		Next, we aim to show that the family $\left\{Y^{\delta}, Z^{\delta}, M^{\delta}\right\}_{\delta \in (0,1]}$ is a Cauchy sequence in $\mathcal{E}^2_{\beta, \mu}$.\\
		Let $0 \leq \delta, \delta' \leq 1$ and define $\widehat{\mathcal{R}} = \mathcal{R}^\delta - \mathcal{R}^{\delta'}$ with $\mathcal{R} \in \left\{Y, Z, M\right\}$.\\
		Using It\^o's formula with expectation, and following the same computations as those used in Proposition \ref{propo2}, along with (\ref{approximating equation}), H\"older's inequality, and the estimations (\ref{estimation for hdelta}) and (\ref{estimation for gdelta}), we can write
		\begin{equation}
				\begin{split}
					&\beta\mathbb{E}\int_{0}^{T}e^{\beta A_{s} + \mu {\kappa}_{s}}|\widehat{Y}_s|^{2} dA_s+\mu\mathbb{E}\int_{0}^{T}e^{\beta A_{s} + \mu {\kappa}_{s}}|\widehat{Y}_s|^{2} d\kappa_s\\
					&+\mathbb{E}\int_{0}^{T}e^{\beta A_{s} + \mu {\kappa}_{s}}\big\|\widehat{Z}_s\big\|^{2}ds+\mathbb{E}\int_{0}^{T} e^{\beta A_{s} + \mu {\kappa}_{s}}d\big[ \widehat{M}\big]_s 		\\
					&\leq 2\left(\delta+\delta'\right)\mathbb{E}\int_{0}^{T}e^{\beta A_{s} + \mu {\kappa}_{s}}\left(\mathfrak{f}_\delta(s,Y^\delta_s) \mathfrak{f}_{\delta'}(s,Y^{\delta'}_s)ds+g_\delta(s,Y^\delta_s) g_{\delta'}(s,Y^{\delta'}_s)d\kappa_s\right)\\
					&\leq 2\left(\delta+\delta'\right)\left( \mathbb{E}\int_{0}^{T}e^{\beta A_{s} + \mu {\kappa}_{s}}\left| \mathfrak{f}_\delta(s,Y^\delta_s)\right|^2 ds \right)^{\frac{1}{2}} \left(  \mathbb{E}\int_{0}^{T}e^{\beta A_{s} + \mu {\kappa}_{s}}\left| \mathfrak{f}_{\delta'}(s,Y^{\delta'}_s)\right|^2 ds\right)^{\frac{1}{2}} \\ &+2\left(\delta+\delta'\right)\left( \mathbb{E}\int_{0}^{T}e^{\beta A_{s} + \mu {\kappa}_{s}}\left| g_\delta(s,Y^\delta_s)\right|^2 d \kappa_s\right)^{\frac{1}{2}} \left(  \mathbb{E}\int_{0}^{T}e^{\beta A_{s} + \mu {\kappa}_{s}} \left| g_{\delta'}(s,Y^{\delta'}_s)\right|^2 d\kappa_s\right)^{\frac{1}{2}} \\
					&\leq 4\left(\delta+\delta'\right)\mathfrak{C}.
				\end{split}
			\label{Cauchy sequnce}
		\end{equation}
		Moreover, following the proof of Proposition \ref{propo1} and using the BDG inequality, we have
		\begin{equation}
			\mathbb{E} \left[\underset{0\leq t\leq T}{\sup}e^{\beta A_{t} + \mu {\kappa}_{t}}|Y_{t}^{\delta}-Y_{t}^{\delta'}|^{2}\right] \leq \mathfrak{C}\left(\delta+\delta'\right).
			\label{Convergence of Y}
		\end{equation}
		Hence, $(Y^\delta, Z^\delta, M^\delta)$ is a Cauchy sequence in $\mathcal{E}^{2}_{\beta, \mu}$. Therefore, there exists a process $(Y, Z, M) \in \mathcal{E}^{2}_{\beta, \mu}$ such that $Y^\delta \to Y$ in $\mathcal{B}^{2}_{\beta, \mu}$, $Z^\delta \to Z$ in $\mathcal{H}^{2}_{\beta, \mu}$, and $M^{\delta} \to M$ in $\mathcal{M}^{2}_{\beta, \mu}$ as $\delta \downarrow 0$.
		
		To complete the proof, it remains to show the convergence of the driver part $\left\{\mathfrak{f}_{\delta}(\cdot, Y^{\delta}_{\cdot}), g_{\delta}(\cdot, Y^{\delta}_{\cdot})\right\}_{\delta \in (0,1]}$ of the GBSDE (\ref{Yosida approximation equation}) to $\left\{\mathfrak{f}(\cdot, Y_{\cdot}), g(\cdot, Y_{\cdot})\right\}$ in an appropriate $\mathbb{L}^2$-space.
		
		From (\ref{estimation for hdelta})–(\ref{estimation for gdelta}), we deduce that $\left\{\delta \, \mathfrak{f}_{\delta}\left(s, Y^{\delta}_s\right) \right\}_{0 < \delta \leq 1}$ belongs to the space $\mathbb{L}^{2}\left(\Omega \times [0,T]; \mathbb{P}(d\omega) \otimes dt \right)$ and that $\left\{\delta \, g_{\delta}\left(s, Y^{\delta}_s\right) \right\}_{0 < \delta \leq 1}$ belongs to the space $\mathbb{L}^{2}\big(\Omega \times [0,T]; \mathbb{P}(d\omega) \otimes d\kappa_t(\omega) \big) := \mathbb{L}^2_{\mathbb{P} \otimes \kappa}$. Here, $\mathbb{P} \otimes d\kappa$ is the positive measure on $\left( \Omega \times [0,T], \mathcal{F} \otimes \mathcal{B}\left([0,T]\right) \right)$ defined for any $\mathcal{V} \in \mathcal{F} \otimes \mathcal{B}\left([0,T]\right)$ by
		$$
		\left( \mathbb{P} \otimes d\kappa \right) \left( \mathcal{V} \right) = \mathbb{E} \left[ \int_{0}^{T} \mathds{1}_{\mathcal{V}}(\omega, s) \, d\kappa_s \right].
		$$
		Furthermore, it holds that $\delta \, \mathfrak{f}_{\delta}\left(s, Y^{\delta}_s\right) \xrightarrow[\delta \to 0^+]{\mathbb{L}^2_{\mathbb{P} \otimes dt}} 0$ and $\delta \, g_{\delta}\left(s, Y^{\delta}_s\right) \xrightarrow[\delta \to 0^+]{\mathbb{L}^2_{\mathbb{P} \otimes d\kappa}} 0$.\\
		By applying the partial reciprocal of the dominated convergence theorem in $\mathbb{L}^2_{\mathbb{P} \otimes dt}$ and $\mathbb{L}^2_{\mathbb{P} \otimes d\kappa}$, respectively, we obtain two subsequences $\left\{\delta_k \, \mathfrak{f}_{\delta_k}\left(s, Y^{\delta_k}_s\right) \right\}_{k \in \mathbb{N}}$ and $\left\{\delta_k \, g_{\delta_k}\left(s, Y^{\delta_k}_s\right) \right\}_{k \in \mathbb{N}}$ of $\left\{\delta \, \mathfrak{f}_{\delta}\left(s, Y^{\delta}_s\right) \right\}_{0 < \delta \leq 1}$ and $\left\{\delta \, g_{\delta}\left(s, Y^{\delta}_s\right) \right\}_{0 < \delta \leq 1}$, respectively, such that
		\begin{equation*}
			\left\lbrace 
			\begin{split}
				&\left\{\delta_k\right\}_{k \in \mathbb{N}} \subset ]0,1]^{\otimes \mathbb{N}},\quad \lim\limits_{k \rightarrow +\infty} \delta_k = 0;\\
				&\delta_k  h_{\delta_k}\big(\omega,s,Y^{\delta_k}_s(\omega)\big) \xrightarrow[k \rightarrow +\infty]{}0,~~\mathbb{P}(d\omega) \otimes dt\mbox{-a.e.;}\\
				&\delta_k g_{\delta_k}\left(\omega,s,Y^{\delta_k}_s(\omega)\right) \xrightarrow[k \rightarrow +\infty]{}0,~~  \mathbb{P}(d\omega) \otimes d\kappa_t(\omega)\mbox{-a.e.}
			\end{split}
			\right. 
		\end{equation*}
		Now, from the convergence (\ref{Convergence of Y}) and the estimation (\ref{Cauchy sequnce}), we can extract a subsequence $\left\{Y^{\delta_k}\right\}_{k \in \mathbb{N}}$ such that $Y^{\delta_k}_t \to Y_t$ a.e. as $k \to +\infty$, with respect to the two positive measures $\mathbb{P} \otimes dt$ and $\mathbb{P} \otimes d\kappa$. Using the continuity of the drivers $\mathfrak{f}$ and $g$ in $y$ (assumption (H7)), we conclude that
		$$
		\lim_{k \to +\infty} \mathfrak{f}_{\delta_k}\left(\omega, t, Y^{\delta_k}_t(\omega)\right) 
		= \lim_{k \to +\infty} \mathfrak{f}(\omega, t, Y^{\delta_k}_t(\omega) + \delta_k h_{\delta_k}) 
		= \mathfrak{f}(\omega, t, Y_t(\omega)),
		$$
		$\mathbb{P}(d\omega) \otimes dt$-a.e., and
		$$
		\lim_{k \to +\infty} g_{\delta_k}\left(\omega, t, Y^{\delta_k}_t(\omega)\right) 
		= \lim_{k \to +\infty} g(\omega, t, Y^{\delta_k}_t(\omega) + \delta_k g_{\delta_k}) 
		= g(\omega, t, Y_t(\omega)),
		$$
		$\mathbb{P}(d\omega) \otimes d\kappa_t(\omega)$-a.e. \\
		Next, using (\ref{estimation for hdelta})–(\ref{estimation for gdelta}) and the Lebesgue dominated convergence theorem, we obtain for any $t \in [0,T]$,
		\begin{equation*}
			\begin{split}
				&\lim\limits_{k \rightarrow +\infty}\Bigg(\mathbb{E}\left[\int_{t}^{T} \left| h_{\delta_k}\left(t,Y^{\delta_k}_s\right) -h(t,Y_s) \right|^2 ds \right]\\
				&\qquad\qquad+\mathbb{E}\left[\int_{t}^{T} \left| g_{\delta_k}\left(s,Y^{\delta_k}_s\right)  -g(s,Y_s) \right|^2 d  \kappa_s \right]\Bigg) =0.
			\end{split}
		\end{equation*}
		Finally, by passing to the limit in the approximating equation (\ref{Yosida approximation equation}) for a subsequence in the corresponding $\mathbb{L}^2$-spaces, we deduce that
		$$
		Y_t=\xi+\int_t^T \mathfrak{f}(s,Y_s)ds+\int_t^T g(s,Y_s)d\kappa_s-\int_t^T Z_sdW_s-\int_t^T dM_s.
		$$
		Completing the proof
	\end{proof}
	\paragraph{General case}
	We may now show the main result of this part after building down the essential basis. Here, the driver $f$ is assumed to be general, i.e., it depends on the $(y,z)$-variables.
	\begin{theorem}\label{thm1}
		Suppose that \textbf{(H-M)$_{\mathbf{2}}$} holds. Then, the GBSDE (\ref{basic GBSDE}) admits a unique $\mathbb{L}^2$-solution.
	\end{theorem}
	\begin{proof}
		The uniqueness is already established in Corollary \ref{coro L2}. To prove existence, we will use a fixed-point argument.\\
		Let us define the space $\mathfrak{L}^2_{\beta,\mu}:=\left(\mathcal{S}^{p,A}_{\beta,\mu}\cap \mathcal{S}^{p,\kappa}_{\beta,\mu}\right) \times \mathcal{H}^p_{\beta,\mu}\times \mathcal{M}^p_{\beta,\mu}$ endowed with the following norm:
		\begin{equation*}
			\begin{split}
				&\Arrowvert (Y,Z,M) \Arrowvert_{\mathfrak{L}^2_{\beta,\mu}}^2 = \mathbb{E}\int_{0}^{T} e^{\beta A_{s} + \mu {\kappa}_{s}}|Y_s|^2 dA_s +\mathbb{E} \int_{0}^{T}e^{\beta A_{s} + \mu {\kappa}_{s}} |Y_s|^2 d\kappa_s\\
				&\hspace*{3.6 cm}+ \mathbb{E}\int_{0}^{T}  e^{\beta A_{s} + \mu {\kappa}_{s}} \left\|Z_s \right\|^2  ds+\mathbb{E}\int_{0}^{T} e^{\beta A_{s} + \mu {\kappa}_{s}} d\left[M\right]_s.
			\end{split}
		\end{equation*}
		It remains to establish the existence, which will be obtained via a fixed point of the contraction mapping $\Psi$ on the space $\left(\mathfrak{L}^2_{\beta, \mu}, \Arrowvert \cdot \Arrowvert_{\mathfrak{L}^2_{\beta, \mu}}\right)$. The function $\Psi$ is defined as follows: for every $(Y, Z, M) \in \mathfrak{L}^2_{\beta, \mu}$, we set $\Psi(Y, Z, M) = (\tilde{Y}, \tilde{Z}, \tilde{M})$, where $(\tilde{Y}, \tilde{Z}, \tilde{M})$ is the unique $\mathbb{L}^2$-solution of the GBSDE (\ref{Special GRBSDE lem}) associated with $(\xi, \mathfrak{f}, g, \kappa)$, where $\mathfrak{f}(t, y) = f(t, y, Z_t)$. The existence of this solution follows from Proposition \ref{propo3}.
		
		Let $(Y,Z,M), (Y',Z',M') \in \mathfrak{L}^2_{\beta,\mu}$ such that $\Psi(Y,Z,M)=(\tilde{Y},\tilde{Z},\tilde{M})$ and $\Psi(Y',Z',M')=(\tilde{Y}',\tilde{Z}',\tilde{M}')$.\\
		Applying It\^o's formula as in Proposition \ref{propo2}, and using assumptions (H2)–(H4) along with Remark \ref{rmq essential}, we obtain, by choosing $\beta \geq 3$ and $\mu \geq 1$,
		\begin{equation*}
			\begin{split}
				&\mathbb{E} \int_{0}^{T} e^{\beta A_{s} + \mu {\kappa}_{s}} \big|\tilde{Y}_s -\tilde{Y}'_s\big|^2 dA_s + \mathbb{E}  \int_{0}^{T}e^{\beta A_{s} + \mu {\kappa}_{s}} \big|\tilde{Y}_s -\tilde{Y}'_s\big|^2 d\kappa_s\\
				&\qquad+ \mathbb{E} \int_{0}^{T} e^{\beta A_{s} + \mu {\kappa}_{s}} \big\|\tilde{Z}_s -\tilde{Z}'_s \big\|^2  ds+\mathbb{E} \int_{0}^{T}e^{\beta A_{s} + \mu {\kappa}_{s}}\big[\tilde{M}-\tilde{M}'\big]_s \\
				&\leq \frac{1}{2} \mathbb{E} \int_{0}^{T}e^{\beta A_{s} + \mu {\kappa}_{s}}\big\|Z_s -Z'_s \big\|^2  ds
			\end{split}
		\end{equation*}
		Then
		\begin{equation*}
			\begin{split}
				&\Arrowvert (\tilde{Y}-\tilde{Y}'),(\tilde{Z}-\tilde{Z}'),(\tilde{M}-\tilde{M}') \Arrowvert_{\mathfrak{L}^2_{\beta,\mu}}^2
				\leq \frac{1}{2} \Arrowvert (Y-Y'),(Z-Z'),(M-M') \Arrowvert_{\mathfrak{L}^2_{\beta,\mu}}^2.
			\end{split}
		\end{equation*}
		Then, $\Psi$ is a contraction mapping on $\left(\mathfrak{L}^2_{\beta, \mu}, \Arrowvert \cdot \Arrowvert_{\mathfrak{L}^2_{\beta, \mu}}\right)$. Hence, there exists a triplet of processes $(Y, Z, M)$ that is a fixed point of $\Psi$, which corresponds to the unique $\mathbb{L}^2$-solution of the GBSDE (\ref{basic GBSDE}).
	\end{proof}

	\subsection{Existence of $\mathbb{L}^p$-solutions for $p \in (1,2)$}
	First, note that the uniqueness result has been established in Corollary \ref{Uniquen}. Using Proposition \ref{propo1}, we obtain an important result that provides an $\mathbb{L}^p$-estimate for the solutions of the GBSDE (\ref{basic GBSDE}) associated with the given data $(\xi, f, g, \kappa)$.
	\begin{corollary}\label{coro Lp}
		Let $(\xi, f, g, \kappa)$ be any set of data satisfying assumption \textbf{(H-M)$_\mathbf{p}$}. For any $\beta, \mu > \frac{2(p-1)}{p}$, there exists a constant $\mathfrak{c}_{\beta, \mu, p, \epsilon}$ such that, whenever $(Y, Z)$ is an $\mathbb{L}^p$-solution of the GBSDE (\ref{basic GBSDE}), we have
		\begin{equation*}
				\begin{split}
					&\mathbb{E}\left[\sup_{t \in [0,T]}e^{\frac{p}{2}\beta A_{t} +\frac{p}{2} \mu {\kappa}_{t}}\big|{Y}_t\big|^p\right]+\mathbb{E} \int_{0}^{T}e^{\frac{p}{2}\beta A_{s} +\frac{p}{2} \mu {\kappa}_{s}}\big|{Y}_s\big|^p dA_s\\
					&+\mathbb{E} \int_{0}^{T}e^{\frac{p}{2}\beta A_{s} +\frac{p}{2} \mu {\kappa}_{s}}\big|{Y}_s\big|^p d\kappa_s+\mathbb{E}\left[\left(\int_{0}^{T}e^{\frac{p}{2}\beta A_{s} +\frac{p}{2} \mu {\kappa}_{s}} \big\|{Z}_s\big\|^2ds\right)^{\frac{p}{2}}\right]\\
					&+\mathbb{E}\left[\left(\int_{0}^{T}e^{\frac{p}{2}\beta A_{s} +\frac{p}{2} \mu {\kappa}_{s}} d\big[{M}\big]_s \right)^{\frac{p}{2}}\right]\\
					\leq & \mathfrak{c}_{\beta,\mu,p,\epsilon}\left( \mathbb{E}\left[ e^{\frac{p}{2}\beta A_{T} +\frac{p}{2} \mu {\kappa}_{T}}\big|{\xi}\big|^p\right] 
					+\int_{0}^{T}e^{\beta A_{s} + \mu {\kappa}_{s}}  \big|\varphi_s\big|^p ds+\mathbb{E}\int_{0}^{T} e^{\beta A_{s} + \mu {\kappa}_{s}}\big|\psi_s\big|^pd\kappa_s\right).
				\end{split}
		\end{equation*}
	\end{corollary}
	
	The main result of this section is stated as follows:
	\begin{theorem}\label{thm2}
		Suppose that \textbf{(H-M)} holds. Then, the GBSDE (\ref{basic GBSDE}) admits a unique $\mathbb{L}^p$-solution.
	\end{theorem}
	\begin{proof}
		We assume that
		\begin{equation}\label{Boudness}
			e^{\frac{p}{2}\beta A_{T} +\frac{p}{2} \mu {\kappa}_{T}}\big|\xi\big|^p+\sup_{t \in [0,T]} e^{\beta A_{t} + \mu {\kappa}_{t}} \big|\varphi_t\big|^p+\sup_{t \in [0,T]} e^{\beta A_{t} + \mu {\kappa}_{t}} \big|\psi_t\big|^p \leq \mathsf{C}.
		\end{equation}
		Note that since $\varphi$ and $\psi$ are $[1, +\infty)$-valued stochastic processes, we deduce that ${\Theta}^{\beta, \mu}_t \leq \mathsf{C}$ for any $t \in [0, T]$. Moreover, from (H5), it follows that 
		$$
		\mathbb{E} \left[ \int_{0}^{T} e^{\beta A_{s} + \mu {\kappa}_{s}} \, (ds + d\kappa_s) \right] < +\infty.
		$$ 
		Additionally, using the definition of the process $\left( {\Theta}^{\beta, \mu}_t \right)_{t \leq T}$, we have ${\Theta}^{\beta, \mu}_t \geq 1$ for any $t \in [0, T]$. Consequently, 
		$$
		\left| \xi \right|^2 + \sup_{t \in [0, T]} \left| \varphi_t \right|^2 + \sup_{t \in [0, T]} \left| \psi_t \right|^2 \leq \mathsf{C}^{\frac{2}{p}}.
		$$
		Therefore, using this and condition (H6), we conclude that
		\begin{equation*}
			\begin{split}
				&\mathbb{E}\left[e^{\beta A_{T} + \mu {\kappa}_{T}}\big|\xi\big|^2+\int_{0}^{T}e^{\beta A_{s} + \mu {\kappa}_{s}}\left|\frac{\varphi_s}{a_s}\right|^2ds+\int_{0}^{T}e^{\beta A_{s} + \mu {\kappa}_{s}}\left|\psi_s\right|^2d\kappa_s\right]\\
				& \leq \mathsf{C}^{\frac{2+p}{p}}+\left(\frac{1}{\epsilon^2}+1\right) \mathsf{C}^{\frac{2}{p}}\mathbb{E}\left[\int_{0}^{T}e^{\beta A_{s} +\mu {\kappa}_{s}}(ds+d\kappa_s)\right]<+\infty.
			\end{split}
		\end{equation*}
		Then, using (\ref{Boudness}) and assumptions (H2)–(H7), we place ourselves within the framework of Theorem \ref{thm1}. Consequently, there exists a unique $\mathbb{L}^2$-solution $(Y^n, Z^n, M^n)$ and, hence, an $\mathbb{L}^p$-solution for any $p \in (1,2)$ (since we are working on the compact interval $[0, T]$) for the GBSDE (\ref{basic GBSDE}). Note also that, by applying Corollary \ref{coro Lp}, we observe that the triplet $(Y, Z, M)$ satisfies the same stated estimation for any $\beta, \mu > \frac{2(p-1)}{p}$.\\
		Using (\ref{Boudness}), we construct a sequence of GBSDEs associated with some data $(\xi_n, f_n, g_n, \kappa)$ such that (\ref{Boudness}) is satisfied and that approximates the GBSDE (\ref{basic GBSDE}). To this end, and to simplify notation, we set $f_0(t) = f(t, 0, 0)$ and $g_0(t) = g(t, 0)$. For each $n \geq 1$, we define a set of data $(\xi_n, f_n, g_n)$ as follows:
		\begin{equation}\label{Def of fn and gn}
			\left\{
			\begin{split}
				\xi_n&=
				\left\{
				\begin{split}
					&\dfrac{\big|\xi\big| \wedge \sqrt[p]{n} 	{\Theta}^{-\frac{\beta}{2} ,-\frac{\mu}{2}}_T}{\big|\xi\big|}\xi &\mbox{ if } \xi \neq 0\\
					&0 &\mbox{ if } \xi = 0
				\end{split}
				\right.\\
				{f}_n(t,y,z)&=
				\left\{
				\begin{split}
					&{f}(t,y,z)-{f}_0(t)+\dfrac{\big|{f}_0(t)\big|\wedge\sqrt[p]{n} {\Theta}^{-\frac{\beta}{p} ,-\frac{\mu}{p}}_t }{\big|{f}_0(t)\big|}{f}_0(t) &\mbox{ if }{f}_0(t) \neq 0\\
					&0 &\mbox{ if } {f}_0(t) = 0
				\end{split}
				\right.\\
				g_n(t,y)&=
				\left\{
				\begin{split}
					&g(t,y)-{g}_0(t)+\dfrac{\big|g_0(t)\big|\wedge\sqrt[p]{n} {\Theta}^{-\frac{\beta}{p} ,-\frac{\mu}{p}}_t }{\big|{g}_0(t)\big|}{g}_0(t) &\mbox{ if } {g}_0(t) \neq 0\\
					&0 &\mbox{ if } {g}_0(t)= 0
				\end{split}
				\right.
			\end{split}
			\right.
		\end{equation}
		
		For each $n \geq 1$, the data $(\xi_n, {f}_n, g_n)$ satisfies condition (\ref{Boudness}). Indeed, it is straightforward to observe the following:
		$$
		e^{\frac{p}{2}\beta A_{T} +\frac{p}{2} \mu {\kappa}_{T}}\big|\xi_n\big|^p+\sup_{t \in [0,T]} e^{\beta A_{t} + \mu {\kappa}_{t}} \big| {f}_n(t,0,0)\big|^p+\sup_{t \in [0,T]} e^{\beta A_{t} + \mu {\kappa}_{t}} \big| {g}_n(t,0)\big|^p \leq n.
		$$
		Therefore, from the previous step, for each $n \geq 1$, there exists a unique triplet $(Y^n, Z^n, M^n)$ that solves the following GBSDE:
		\begin{equation}\label{basic GBSDE--appr}
			Y^n_t=\xi_n+\int_{t}^{T}{f}_n(s,Y^n_s,Z^n_s)ds+\int_{t}^{T}g_n(s,Y^n_s)d\kappa_s-\int_{t}^T Z^n_s d W_s-\int_{t}^T dM^n_s
		\end{equation}	
		for any $t \in [0,T]$.
		
		Let $n \geq m  \geq 1$. Set $\widehat{R}=\mathcal{R}^n-\mathcal{R}^m$ for $\mathcal{R} \in \left\{\xi, Y, Z, M\right\}$. Using again Corollary \ref{used in a priori estimates}, we get
		\begin{equation}\label{following}
				\begin{split}
					&e^{\frac{p}{2}\beta A_{t} +\frac{p}{2} \mu {\kappa}_{t}}\big|\widehat{Y}_t\big|^p+\frac{p}{2}\beta \int_{t}^{T}e^{\frac{p}{2}\beta A_{s} +\frac{p}{2} \mu {\kappa}_{s}}\big|\widehat{Y}_s\big|^pdA_s+\frac{p}{2}\mu  \int_{t}^{T}e^{\frac{p}{2}\beta A_{s} +\frac{p}{2} \mu {\kappa}_{s}}\big|\widehat{Y}_s\big|^pd\kappa_s\\
					&+c(p)\int_{0}^{t}e^{\frac{p}{2}\beta A_{s} +\frac{p}{2} \mu {\kappa}_{s}}\big|\widehat{Y}_s\big|^{p-2}\mathds{1}_{Y_s \neq 0}\big\|\widehat{Z}_s\big\|^2ds\\
					&+c(p)\int_{t }^{T}e^{\frac{p}{2}\beta A_{s} +\frac{p}{2} \mu {\kappa}_{s}}\big|\widehat{Y}_s\big|^{p-2}\mathds{1}_{Y_s \neq 0}d \big[\widehat{M}\big]^c_s\\
					\leq &e^{\frac{p}{2}\beta A_{T} +\frac{p}{2} \mu {\kappa}_{T}}\big|\widehat{\xi}\big|^p
					+p\int_{t}^{T}e^{\frac{p}{2}\beta A_{s} +\frac{p}{2} \mu {\kappa}_{s}} \big|\widehat{Y}_s\big|^{p-1} \check{\widehat{Y}}_s \big(f_n(s,Y^n_s,Z^n_s)-f_m(s,Y^m_s,Z^m_s)\big) ds\\
					&+p\int_{t}^{T} e^{\frac{p}{2}\beta A_{s} +\frac{p}{2} \mu {\kappa}_{s}}\big|\widehat{Y}_s\big|^{p-1} \check{\widehat{Y}}_s \big(g_n(s,Y^n_s)-g_m(s,Y^m_s)\big)d\kappa_s\\
					&-p\int_{t}^{T}e^{\frac{p}{2}\beta A_{s} +\frac{p}{2} \mu {\kappa}_{s}} \big|\widehat{Y}_s\big|^{p-1} \check{\widehat{Y}}_s \widehat{Z}_s dW_s-p\int_{t }^{T}{\Theta}^{\frac{p}{2}\beta ,\frac{p}{2}\mu}_s \big|\widehat{Y}_{s-}\big|^{p-1} \check{\widehat{Y}}_{s-} d\widehat{M}_s\\
					&-\sum_{t  <s \leq T} e^{\frac{p}{2}\beta A_{s} +\frac{p}{2} \mu {\kappa}_{s}}\left\{\big|\widehat{Y}_{s-}+\Delta \widehat{M}_s\big|^p-\big|\widehat{Y}_{s-}\big|^p-p\big|\widehat{Y}_{s-}\big|^{p-1}\check{\widehat{Y}}_{s-}\Delta \widehat{M}_s\right\}.
				\end{split}
		\end{equation}
		From (\ref{Def of fn and gn}) and assumptions (H2)-(H4) on $f$ and $g$, we have
		\begin{equation*}
			\begin{split}
				&\widehat{Y}_t \left({f}_n(s,Y^n_s,Z^n_s)-{f}_m(s,Y^m_s,Z^m_s)\right)\\
				&=\widehat{Y}_t \left({f}(s,Y^n_s,Z^n_s)-{f}(s,Y^m_s,Z^m_s)\right)+\widehat{Y}_t \left({f}_n(s,0,0)-{f}_m(s,0,0)\right)\\
				& \leq  \frac{c(p)}{2}\|\widehat{Z}_s\|^2+\big|\widehat{Y}_t\big|\big|{f}_n(s,0,0)-{f}_m(s,0,0)\big|
			\end{split}
		\end{equation*}
		and
		\begin{equation*}
			\begin{split}
				\widehat{Y}_t \left({g}_n(s,Y^n_s)-{g}_m(s,Y^m_s)\right)
				&=\widehat{Y}_t \left({g}(s,Y^n_s)-{g}(s,Y^m_s)\right)+\widehat{Y}_t \left({g}_n(s,0)-{g}_m(s,0)\right)\\
				& \leq  \big|\widehat{Y}_s\big|\big|{g}_n(s,0)-{g}_m(s,0)\big|
			\end{split}
		\end{equation*}
		Then, we obtain an analogous estimation to (\ref{CV}) for the driver $f_n$ and a simpler one than (\ref{CV1}) for the coefficient $g$. Following this, using (\ref{following}) and re-performing the calculations from Proposition \ref{propo1}, we deduce that, for any $\beta, \mu > \frac{2(p-1)}{p}$, there exists a constant $\mathfrak{c}_{\beta, \mu, p, \epsilon}$ such that
		\begin{equation}\label{RHS n}
			\begin{split}
				&\mathbb{E}\left[\sup_{t \in [0,T]}e^{\frac{p}{2}\beta A_{t} +\frac{p}{2} \mu {\kappa}_{t}}\big|\widehat{Y}_t\big|^p\right]+\mathbb{E} \int_{0}^{T}e^{\frac{p}{2}\beta A_{s} +\frac{p}{2} \mu {\kappa}_{s}}\big|\widehat{Y}_s\big|^pdA_s\\
				&\quad+\mathbb{E} \int_{0}^{T}e^{\frac{p}{2}\beta A_{s} +\frac{p}{2} \mu {\kappa}_{s}}\big|\widehat{Y}_s\big|^pd\kappa_s+\mathbb{E}\left[\left( \int_{0}^{T}e^{\beta A_{s} + \mu {\kappa}_{s}}\big\|\widehat{Z}_s\big\|^2ds\right)^{\frac{p}{2}}\right]\\
				&+\mathbb{E}\left[ \left( \int_{0}^{T}e^{\beta A_{s} + \mu {\kappa}_{s}} d\big[\widehat{M}\big]_s\right)^{\frac{p}{2}}\right]  \\
				\leq & \mathfrak{c}_{\beta,\mu,p,\epsilon}\left(\mathbb{E}\left[e^{\frac{p}{2}\beta A_{T} +\frac{p}{2} \mu {\kappa}_{T}} \big|\widehat{\xi}\big|^p\right]  +\mathbb{E} \int_{0}^{T}e^{\beta A_{t} + \mu {\kappa}_{s}}  \big|{f}_n(s,0,0)-{f}_m(s,0,0)\big|^pds \right.\\
				&\left. \qquad\qquad+\mathbb{E} \int_{0}^{T}e^{\beta A_{s} +\mu {\kappa}_{s}} \big|{g}_n(s,0)-{g}_m(s,0)\big|^pd\kappa_s\right).
			\end{split}
		\end{equation}
		By using the basic inequality 
		\begin{equation*}
			\left(\sum_{i=1}^n |X_i|\right)^p \leq n^p \sum_{i=1}^n |X_i|^p \qquad \forall (n, p) \in \mathbb{N}^\ast \times (0, +\infty),
		\end{equation*}
		along with (\ref{Def of fn and gn}) and assumption (H5), we obtain
		\begin{equation}\label{with Coro}
			\left\{
			\begin{split}
				\big|\widehat{\xi}\big|^p &\leq 2^{p+1} \big|{\xi}\big|^p, \\
				\big|{f}_n(s, 0, 0) - {f}_m(s, 0, 0)\big|^p &\leq 2^{p+1} \left|\varphi_s\right|^p, \quad \mathbb{P} \otimes dt \text{-a.e.}, \\
				\big|{g}_n(s, 0) - {g}_m(s, 0)\big|^p &\leq 2^{p+1} \left|\psi_s\right|^p, \quad \mathbb{P} \otimes d\kappa_t \text{-a.e.}
			\end{split}
			\right.
		\end{equation}
		Since $\lim\limits_{n \to +\infty} f_n(t, 0, 0) = f_0(t)$ $\mathbb{P} \otimes dt$-a.e. and $\lim\limits_{n \to +\infty} g_n(t, 0) = g_0(t)$ $\mathbb{P} \otimes d\kappa_t$-a.e., it follows from \textsc{(H1)} and \textsc{(H5)} that we can apply the Lebesgue dominated convergence theorem. Hence, we deduce that the right-hand side of (\ref{RHS n}) tends to zero as $n, m \to +\infty$. Therefore, the left-hand side of (\ref{RHS n}) also tends to zero. Consequently, we derive the following convergences:
		\begin{equation}\label{CVV}
			\begin{split}
				\lim\limits_{n,m \rightarrow +\infty}\left(\big\|{Y}^n-Y^m\big\|^p_{\mathfrak{B}^p_{\beta,\mu}}+\big\|{Z}^n-Z^m\big\|^p_{\mathcal{H}^p_{\beta,\mu}}+\big\|{M}^n-M^m\big\|^p_{\mathcal{M}^p_{\beta,\mu}}\right)
				=0.
			\end{split}
		\end{equation}
		Hence, $\left\{(Y^n, Z^n, M^n)\right\}_{n \geq 1}$ is a Cauchy sequence in the Banach space $\mathcal{E}^p_{\beta, \mu}$ for any $\beta, \mu > \frac{2(p-1)}{p}$. It then converges to a process $(Y, Z, M) \in \mathcal{E}^p_{\beta, \mu}$. Moreover, using (\ref{with Coro}) and Corollary \ref{coro Lp}, we deduce that, for any $\beta, \mu > \frac{2(p-1)}{p}$, there exists a constant $\mathfrak{c}_{\beta, \mu, p, \epsilon}$ (independent of $n$) such that
		\begin{equation}\label{UI}
				\begin{split}
					&\mathbb{E}\left[\sup_{t \in [0,T]}e^{\frac{p}{2}\beta A_{t} +\frac{p}{2} \mu {\kappa}_{t}}\big|{Y}^n_t\big|^p\right]+\mathbb{E} \int_{0}^{T}e^{\frac{p}{2}\beta A_{s} +\frac{p}{2} \mu {\kappa}_{s}}\big|{Y}^n_s\big|^pdA_s\\
					&+\mathbb{E} \int_{0}^{T}e^{\frac{p}{2}\beta A_{s} +\frac{p}{2} \mu {\kappa}_{s}} \big|{Y}^n_s\big|^pd\kappa_s
					+\mathbb{E}\left[\left( \int_{0}^{T}e^{\beta A_{s} +\mu {\kappa}_{s}}\big\|{Z}^n_s\big\|^2ds\right)^{\frac{p}{2}}\right]\\
					&+\mathbb{E}\left[ \left( \int_{0}^{T}e^{\beta A_{s} + \mu {\kappa}_{s}} d\big[M^n\big]_s\right)^{\frac{p}{2}}\right] \\
					\leq & \mathfrak{c}_{\beta,\mu,p,\epsilon}\left(\mathbb{E}\left[e^{\frac{p}{2}\beta A_{T} +\frac{p}{2} \mu {\kappa}_{T}}\big|{\xi}\big|^p\right]  +\mathbb{E} \int_{0}^{T}e^{\beta A_{s} +\mu {\kappa}_{s}}  \big|\varphi_s\big|^pds +\mathbb{E} \int_{0}^{T}e^{\beta A_{s} + \mu {\kappa}_{s}} \big|\psi_s\big|^pd\kappa_s\right).
				\end{split}
		\end{equation}
		
		It remains to confirm that the limiting process solves the generalized BSDE (\ref{basic GBSDE}). To this end, since $p > 1$, we apply the BDG inequality along with (\ref{CV}) to obtain
		\begin{equation*}
			\begin{split}
				&\mathbb{E}\left[\sup_{t \in [0,T]} \left|\int_{t}^{T}Z^n dW_s-\int_{t}^{T}Z_s dW_s\right|^p\right]\\ 
				&\leq \mathfrak{c}\mathbb{E}\left[\left(\int_{0}^{T}\big\|Z^n_s-Z_s\big\|^2 ds\right)^{\frac{p}{2}}\right]\xrightarrow[n \rightarrow +\infty]{}0
			\end{split}
		\end{equation*}
		and
		\begin{equation*}
			\begin{split}
				&\mathbb{E}\left[\sup_{t \in [0,T]} \left|\int_{t}^{T}dM^n_s -\int_{t}^{T} dM_s\right|^p\right]
				\leq \mathfrak{c}\mathbb{E}\left[\left(\int_{0}^{T}d\big[M^n-M\big]_s \right)^{\frac{p}{2}}\right]\xrightarrow[n \rightarrow +\infty]{}0.
			\end{split}
		\end{equation*}
		From assumptions (H4) and (H7), we have
		\begin{equation}\label{semiCV}
			\begin{split}
				\left| f(t, Y^n_t, Z^n_t) - f(t, Y_t, Z_t) \right| 
				&\leq \eta_t \|Z^n_t - Z_t\| + \left| f(t, Y^n_t, Z_t) - f(t, Y_t, Z_t) \right| \\
				&~~~\xrightarrow[n \to +\infty]{} 0.
			\end{split}
		\end{equation}
		Moreover, from assumptions (H5)–(H6) and Jensen's inequality, we have
		\begin{equation*}
			\begin{split}
				&\mathbb{E}\int_{0}^{T}\left|\frac{f(t,Y^n_t,Z_t)-f(t,Y_t,Z_t) }{a_s}\right|^p ds\\
				& \leq 2^p \left(\mathbb{E}\int_{0}^{T}\left|\frac{f(t,Y^n_t,Z_t) }{a_s}\right|^p ds+\mathbb{E}\int_{0}^{T}\left|\frac{f(t,Y_t,Z_t) }{a_s}\right|^p ds\right)\\
				& \leq 2^{2p} \left(\mathbb{E}\int_{0}^{T}\left|\frac{f(t,Y^n_t,Z_t)-f(t,Y^n_t,0) }{a_s}\right|^p ds+\mathbb{E}\int_{0}^{T}\left|\frac{f(t,Y^n_t,0) }{a_s}\right|^p ds \right.\\
				& \left.\qquad\qquad+\mathbb{E}\int_{0}^{T}\left|\frac{f(t,Y_t,Z_t)-f(t,Y_t,0)}{a_s}\right|^p ds +\mathbb{E}\int_{0}^{T}\left|\frac{f(t,Y_t,0) }{a_s}\right|^p ds
				\right)\\
				& \leq 2^{2p}\left(2\mathbb{E}\int_{0}^{T}\left\|Z_s\right\|^pds+\frac{2}{\epsilon^p}\mathbb{E}\int_{0}^{T}\left|\varphi_s\right|^p ds 
				+\mathbb{E}\int_{0}^{T}\left(\left|Y^n_s\right|^p  +\left|Y_s\right|^p\right)  ds  \right)\\
				& \leq 2^{2p}\left(2T^{\frac{2-p}{2}}\mathbb{E}\left[\left(\int_{0}^{T}\left\|Z_s\right\|^2ds\right)^{\frac{p}{2}}\right]+\frac{2}{\epsilon^p}\mathbb{E}\int_{0}^{T}\left|\varphi_s\right|^p ds \right.\\
				&\left.\qquad\qquad+\frac{2}{\epsilon^{2}}\mathbb{E}\int_{0}^{T}\left(\left|Y^n_s\right|^p +\left|Y_s\right|^p\right)  dA_s \right).
			\end{split}
		\end{equation*}
		Using (\ref{Def of fn and gn}), (\ref{CVV}), (\ref{UI}), (\ref{semiCV}), and the Lebesgue dominated convergence theorem, we obtain
		$$
		\lim\limits_{n \rightarrow +\infty}\mathbb{E}\int_{0}^{T}\left|\frac{f_n(t,Y^n_t,Z^n_t)-f(t,Y_t,Z_t) }{a_s}\right|^p ds=0.
		$$
		Similarly, we can show that
		\begin{equation*}
			\begin{split}
				&\lim\limits_{n \rightarrow +\infty}\mathbb{E}\int_{0}^{T}\left|g_n(s,Y^n_s)-g(s,Y_s)\right|^p d\kappa_s=0.
			\end{split}
		\end{equation*}
		Finally, by passing to the limit term by term in (\ref{basic GBSDE--appr}), we deduce that the limiting process $(Y, Z, M)$ is the $\mathbb{L}^p$-solution of the GBSDE (\ref{basic GBSDE}).
		
		This completes the proof.
	\end{proof}

	

	\bibliography{Ref}
	\bibliographystyle{abbrv}
\end{document}